\documentclass[twocolumn]{autart}
\usepackage{amsmath}
\usepackage{amssymb}

\usepackage{epsfig,amssymb,latexsym,amsmath}
\usepackage{euscript}
\usepackage{upgreek}
\usepackage{color}

\usepackage[latin1]{inputenc}
\usepackage[T1]{fontenc}
\usepackage[english]{babel}
\usepackage{dsfont}
\usepackage{fancyhdr}
\usepackage{ amsfonts, amssymb}

\newtheorem{assumption}{Assumption}
\newtheorem{proposition}{Proposition}
\newtheorem{remark}{Remark}
\newtheorem{theorem}{Theorem}
\newtheorem{definition}{Definition}
\newtheorem{lemma}{Lemma}
\newtheorem{corollary}{Corollary}

\newcommand{\vph}{\varphi}

\newcommand{\s}{\sigma}

\newcommand{\intt}{\operatorname{int\,}}

\newcommand{\T}{\operatorname{T\,}}

\newcommand{\id}{\operatorname{id\,}}

\newcommand{\diag}{\operatorname{diag\,}}

\newcommand{\loc}{\operatorname{loc\,}}

\newcommand{\wt}{\widetilde}

\begin{document}

\begin{frontmatter}
	
	\title{Dwell-time stability conditions for infinite dimensional impulsive systems}
	
	\thanks[footnoteinfo]{This work was supported by the German Research Foundation (DFG), Grant number DA767/12-1. Corresponding author: Sergey Dashkovskiy.}
	
	\author[First]{Sergey Dashkovskiy}
	\author[First]{Vitalii Slynko}

	\address[First]{Institute of Mathematics, University of Wuerzburg, Emil-Fischer-Str. 40, 97074 Wuerzburg, Germany\\ email: \{sergey.dashkovskiy,vitalii.slynko\}@mathematik.uni-wuerzburg.de}
	
	\begin{abstract}                
We consider nonlinear impulsive systems on Banach spaces subjected to disturbances and look for dwell-time conditions guaranteeing the the ISS property. In contrary to many existing results our conditions cover the case where both continuous and discrete dynamics can be unstable simultaneously. Lyapunov type methods are use for this purpose. The effectiveness of our approach is illustrated on a rather nontrivial example, which is feedback connection of an ODE and a PDE systems. 
	\end{abstract}

\end{frontmatter}
\section{Introduction}
Hybrid systems accommodate continuous and discontinuous behavior, which allow to model modern practical processes, where a combination of analog and digital effects takes place as well as other processes where, for example, collisions can change the systems state instantaneously. The general theory of hybrid systems includes results on existence, uniqueness, continuous dependence on initial data, stability and robustness of solutions, see \cite{michel,teel,liberzon}. Impulsive systems are a particular subclass of hybrid ones, see \cite{sam-per,lak-bainov}, where many stability results were developed, in particular by means of the linear approximations and by the direct Lyapunov method as well as by means of the corresponding comparison principle. Such systems can be seen as a combination of continuous and discrete subsystems. Stability conditions for nonlinear impulsive systems using continuously differential Lyapunov functions were provided in  \cite{sam-per}. The proofs are based on the comparison principle and require that either of the continuous and discrete subsystems is stable. More general stability conditions were developed in  \cite{lak-bainov} on the base of Lyapunov functions discontinuous in time. These conditions allow also to establish stability even in the case, when both dynamics are unstable. Moreover in this class of Lyapunov functions these conditions are necessary and sufficient \cite{ignat-NA}. Practical applications of these results lead to difficulties in contructions of suitable Lyapunov functions. This problem was considered in \cite{Sl-PM,briat} for the case of linear finite dimensional impulsive systems with constant coefficients. Even in this relatively simple case these works need to use lyapunov functions, which depend explicitly on time. In general, a construction of Lyapunov function leads to a rather complex boundary value problem for a systems of matrix differential equations \cite{briat}. 

In contrary to the above mentioned  works, we will consider nonlinear infinite dimensional impulsive systems. Our aim is to establish stability conditions by means of smooth Lyapunov functions, which do not depend on time at least in case when the right hand sides of the equations do not depend on time. As well we would like to cover the case, where both continuous and discrete dynamics can be unstable. Our class of systems is essentially more general, than in the above mentioned literature and the existing methods cannot be applied directly for our purposes. 

Certain results in this direction exist: The work \cite{liu} uses second order and the work \cite{Dv-Sl-Nosc} uses even higher order derivatives of Lyapunov functions, which leads to additional restrictions on the equations constituting the impulsive system. Different related approaches were developed in 
\cite{SBT-SCL,SB-MS}, devoted to linear impulsive systems on Banach spaces containing continuous operators only. Also these results cannot be extended directly to the class of systems considered in this paper.

Furthermore, we will consider systems with disturbances and hence we will work in the ISS framework. This allows for potential applications in view of studying stability of interconnected systems. Several related results in this framework exist, see \cite{dash-miron13,miron-liberzon18,dash-feketa}, for example.
The work \cite{dash-miron13} provides stability conditions of the dwell-time type for systems in Banach spaces, however it was assumed that either discrete or continuous dynamics possesses the ISS property. Interconnections of impulsive systems, such that not all of them are ISS were studied in \cite{miron-liberzon18}. The work \cite{dash-feketa}  studies ISS properties of impulsive systems where jump mappings may depend on time and provides dwell-time conditions establishing the ISS property. However only finite dimensional systems were considered in the last two papers.

The main result of our work extends the results of \cite{dash-miron13} to the general case, where we do not require that either continuous or discrete dynamics is ISS. To this end we use completely different approach. We split the state space into two subsets $X=G^-\cup G^+$ with $\text{int }G^-\cap G^+=\emptyset$ and such that $G^+$ is forward invariant. We use two auxiliary functions, one is $V$ of the Lyapunov type the other one is $W$ of Chetaev type.  The latter one is used to establish the invariance property of $G^+$. 

The Lyapunov function $V$ decays along the trajectories in $G^-$ but can increase in $G^+$ and allows to estimate the change of solutions between the jumps. It is assumed that the jumps from the state in $G^+$ are always stabilizing, and there is no such restriction in $G^-$. Finally, the ISS property is guaranteed by a dwell-time condition. In the particular case, when both continuous and discrete dynamics are not stable this conditions restricts the jumps frequency both from below and above.

We apply our result to a feedback connection of a linear ODE and a nonlinear PDE of the parabolic type. This illustrates how our approach can be applied and demonstrates its powerfulness. Let us note that this example cannot be handled in view of stability by the existing results because of its nonlinearity, possible instability of both dynamics as well as irregular time instants of jumps. 
\section{Preliminaries}
We will use the following classes of continuous functions, frequently called comparison functions \\
$\mathcal P=\{\gamma:\,\Bbb R_+\to\Bbb R_+|\,\gamma(r)=0\;\Leftrightarrow\; r=0\},$\\
$\mathcal K=\{\gamma\in\mathcal P|\,\gamma\text{ is strictly increasing}\},$\\
$\mathcal K_{\infty}=\{\gamma\in\mathcal K|\,\gamma\text{ is unbounded}\},$\\
$\mathcal K_{\infty}^2=\{\gamma\,:\,\Bbb R_+^2\to\Bbb R_+\,|\,\gamma(\cdot,s)\in\mathcal K_{\infty},\,\gamma(s,\cdot)\in\mathcal K_{\infty}\},$\\
$\mathcal L=\{\gamma:\Bbb R_+\to\Bbb R_+|\,\text{decreasing with}\lim\limits_{t\to\infty}\gamma(t)=0\},$\\
$\mathcal K\mathcal L=\{\beta:\Bbb R_+\times\Bbb R_+\to\Bbb R_+|\,\,\beta(\cdot,t)\in\mathcal K,\,\beta(s,\cdot)\in\mathcal L\}.$
For any $\varrho\in\mathcal K$ and any $a,b\in\Bbb R_+$ it follows by the monotonicity that
\begin{equation}\label{CF}
\varrho(a+b)\le\varrho(2a)+\varrho(2b).
\end{equation}
$C[0,l]$ denotes the space of functions continuous on $[0,l]$ with values in $\Bbb R$ and the norm
$\|f\|_{C[0,l]}=\max\limits_{x\in [0,l]}|f(x)|$, \\
$C^k[0,l]$ stands for the space of $k$-times continuously differentiable functions normed by
$\|f\|_{C^k[0,l]}=\max\limits_{p=0,\dots,k}\max\limits_{x\in [0,l]}|f^{(p)}(x)|$.
$H^0[0,l]=L^2[0,l]$ is the Hilbert space of measurable and square integrable functions with the scalar product $(f,g)_{L^2[0,l]}=\int\limits_0^lf(z)g(z)\,dz$. $\mathfrak L(L^2[0,l])$ denotes the Banach algebra of linear bounded operators on $L^2[0,l]$.

The Hilbert space $H^k[0,l]\subset L^2[0,l]$ is a subset of $L^2[0,l]$ of functions $f$ such that  $f^{(p)}\in L^2[0,l]$ for , $p=0,\dots,k$ and the scalar product defined by
\begin{equation}\label{0.1}
	(f,g)_{H^k[0,l]}=\sum\limits_{p=0}^k\int\limits_0^lf^{(p)}(z)g^{(p)}(z)\,dz.
\end{equation}
Recall that $H^k[0,l]$ is the completion of $C^k[0,l]$ with respect to the norm  $\|f\|_{H^k[0,l]}=\sqrt{(f,f)_{H^k[0,l]}}$.
$C_0^{\infty}[0,l]$ is the space of infinitely smooth on $[0,l]$ functions vanishing in the vicinity of
$x=0$ and $x=l$. $C_0^{\infty}([0,T],C_0^{\infty}[0,l])$ is the set of mappings  $f\,:\,[0,T]\to C_0^{\infty}[0,l]$ vanishing in the vicinity of  $t=0$ and $t=T$. Completion of $C_0^{\infty}[0,l]$ with respect to the norm  $\eqref{0.1}$ is denoted by $H_0^k[0,l]$. In the space $H_0^1[0,l]$, due to the inequality of Friedrichs, the norms
$\|\cdot\|_{H^1[0,l]}$ and $\|\cdot\|_{H^1_0[0,l]}$ defined by
$\|f\|_{H^1_0[0,l]}=\|f_z\|^2_{L^2[0,l]}$, are equivalent.
$L^{\infty}(\Bbb R_+)$ is the space of measurable and essentially bounded functions $f:\Bbb R_+\to\Bbb R$,

For $M\subset\Bbb R$ and a Banach space $X$ by $L^{\infty}(M,X)$ we denote the space of mappings $f\,:\,M\to X$ normed by $\|f\|_{L^{\infty}}=\sup_{m\in M}\|f(m)\|_{X}$, where in the particular case
$M=\Bbb Z_+$ we write $L^{\infty}(\Bbb Z_+,X)=l^{\infty}(X)$. $B_r(x)$ denotes the open ball centered at
$x\in X$ of radius $r>0$. For Banach spaces $X, Y$ and mappings $f\,:\,\Bbb R\to X$, $g\,:\,\Bbb R\to Y$ the Cartesian product  $f\times g\,:\,\Bbb R\to X\times Y$ is defined by $(f\times g)(t)=(f(t),g(t))$. For $\alpha\in(0,1]$  the space of locally H\"older continuous functions $=f:\Bbb R_+\to X$ is denoted by $H_\alpha(\Bbb R_+,X)$ and let $H_{loc}(\Bbb R_+,X):=\cup_{\alpha\in(0,1]}H_\alpha(\Bbb R_+,X)$.

For a linear bounded operator $A$ defined on a Bacnach space by $\s(A)$ we denote its spectrum and by
$r_{\s}(A)$ its spectral radius.

$\Bbb R^{n\times m}$ denotes the linear space of matrices of the size
$n\times m$, where in case $m=n$, the space $\Bbb R^{n\times n}$ is a Banach algebra.
For $A\in\Bbb R^{n\times n}$ we denote $\Bbb R^n$: $\|A\|=\sup_{\|x\|=1}\|Ax\|=\lambda_{\max}^{1/2}(A^{\T}A)$.
$\Bbb S^n$ denotes the set of symmetric matreces of the size $n$. For $P,Q\in\Bbb S^n$ we write $P\succ Q$ if the matrix $P-Q$ is positive definite.
For $A\in\Bbb S^{n}$ by $\lambda_{\min}(A)$ and $\lambda_{\max}(A)$ we denote the smallest and the largest eigenvalues of $A$ respectively. The inequality $P\le Q$ should be understood element-wise.

The following well-known inequalities will be used
\begin{equation}\label{Yang}
\gathered
xy\le \frac{x^{p_1}}{p_1}+\frac{y^{p_2}}{p_2},\quad x\ge 0,\quad y\ge 0,\quad p_1\in(0,\infty),\quad \\
\frac{1}{p_1}+\frac{1}{p_2}=1,
\endgathered
\end{equation}
and for any $f\in L^{p_1}[0,l]$ and $g\in L^{p_2}[0,l]$ it holds that
\begin{equation}\label{Holder}
\gathered
\Big|\int\limits_{0}^lf(t)g(t)\,dt\Big|\le\Big(\int\limits_0^l|f(t)|^{p_1}\,dt\Big)^{1/p_1}\Big(\int\limits_0^l|g(t)|^{p_2}\,dt\Big)^{1/p_2},\\
\quad p_1\in(1,\infty),\quad \frac{1}{p_1}+\frac{1}{p_2}=1,
\endgathered
\end{equation}
known as Young's and H\"older's inequalities, in the particular case $p_1=p_2=2$ they are also known as Cauchy inequalities.
For any $f\in H_0^1(0,l)$ also holds
\begin{equation}\label{Fr1}
\gathered
\|f\|_{H^1_0(0,l)}^2\ge\frac{\pi^2}{l^2}\|f\|_{L^2[0,l]}^2,
\endgathered
\end{equation}
and if additionally $f\in H_0^1(0,l)\cap H^2(0,l)$, then
\begin{equation}\label{Fr2}
\gathered
\|\partial_{zz}f\|_{L^2[0,l]}^2\ge\frac{\pi^2}{l^2}\|\partial_z f\|_{L^2[0,l]}^2,
\endgathered
\end{equation}
\section{Stability notions}
Let us introduce dynamical systems that we will consider adapting definitions from \cite{sontag,MKK19,MiP20}.
\begin{definition}
	Let $X$ be the state space with the norm $\|\cdot\|_X$ and
	$\mathcal U_1\subset\{f\,:\Bbb R\to U_1\}$ be the space of input signals normed by  $\|\cdot\|_{\mathcal U_1}$
	with values in a nonempty subset $U_1$ of some linear normed space and invariant under the time shifts, that is if $d_1\in\mathcal U_1$ and $\tau\in\Bbb R$, then $\mathcal S_{\tau}d_1\in\mathcal U_1$, where $\mathcal S_s\,:\,\mathcal U_1\to \mathcal U_1$, $s\in\Bbb R$ is the linear operator defined by $\mathcal S_s u(t)=u(t+s)$ and satisfying $\|\mathcal S_s\|\le1$. It is also assume that for all $u_1,u_2\in\mathcal U$ and any $t\ge t_0$ we have that
	$u(\tau):=u_1(\tau)$ for $\tau\in[t_0,t]$ and $u(\tau):=u_2(\tau)$ for $\tau>t_0$ it holds that $u\in\mathcal U$.
	
	The triple $\Sigma_c=(X,\mathcal U_1,\phi_c)$ is called dynamical system with inputs if the mapping
	$\phi_c\,\,:\,(t,t_0,x,d_1)\mapsto\phi_c(t,t_0,x,d_1)$ defined for all $(t,t_0,x,d_1)\in[t_0,t_0+\epsilon_{t_0,x,d_1})\times\Bbb R\times X\times\mathcal U_1$ for some positive  $\epsilon_{t_0,x,d_1}$ and satisfies the following axioms
	
	($\Sigma_c1$) for $t_0\in\Bbb R$, $x\in X$, $d_1\in\mathcal U_1$, $t\in [t_0,t_0+\epsilon_{t_0,x,d_1})$ the value of
	$\phi_c(t,t_0,x,d_1)$ is well defined and the mapping $t\mapsto\phi_c(t,t_0,x,d_1)$ is continuous on  $(t_0,t_0+\epsilon_{t_0,x,d_1})$ with $\lim_{t\to t_0+}\phi_c(t,t_0,x,d_1)=x$;
	
	($\Sigma_c2$) $\phi_c(t,t,x,d_1)=x$ for any $(x,d_1)\in X\times\mathcal U_1$, $t\in\Bbb R$
	
	($\Sigma_c3$) for any $t_0\in\Bbb R$, $(t,x,d_1)\in[t_0,t_0+\epsilon_{t_0,x,d_1})\times X\times\mathcal U_1$
	and $\wt d_1\in\mathcal U_1$ with $d_1(s)=\wt d_1(s)$ for $s\in[t_0,t]$ it holds that $\phi_c(t,t_0,x,d_1)=\phi_c(t,t_0,x,\wt d_1)$;
	
	($\Sigma_c4$) for any $(x,d_1)\in X\times\mathcal U_1$ and $t\ge\tau\ge t_0$ with $\tau\in[t_0,t_0+\epsilon_{t_0,x,d_1})$, $t\in[\tau,\tau+\epsilon_{\tau,\phi(\tau,t_0,x,d_1),d_1})\cap[t_0,t_0+\epsilon_{t_0,x,d_1})$
	it holds that
	\begin{equation*}
	\phi_c(t,t_0,x,d_1)=\phi_c(t,\tau,\phi(\tau,t_0,x,d_1),d_1),
	\end{equation*}
		($\Sigma_c5$) for any $(x,d_1)\in X\times\mathcal U_1$ and $t\in [t_0,t_0+\epsilon_{t_0,x,d_1})$,
	it holds that
	\begin{equation*}
	\gathered
	\epsilon_{t_0+\tau,x,d_1}=\epsilon_{t_0,x,\mathcal S_{\tau}d_1},\\
	\phi_c(t+\tau,t_0+\tau,x,d_1)=\phi_c(t,t_0,x,\mathcal S_{\tau}d_1).
	\endgathered
	\end{equation*}
\end{definition}
Note that $(\Sigma_c5)$ implies that for all $t\in[\tau,\tau+\epsilon_{\tau,x_0,d_1})$, $\tau\le t$
\begin{equation}\label{GrPr}
\phi_c(t,\tau,x,d_1)=\phi_c(t-\tau,0,x,\mathcal S_{\tau}d_1).
\end{equation}
Systems with impulsive actions are defined as follows
\begin{definition}
	Let $\mathcal E=\{\tau_k\}_{k=0}^{\infty}, \tau_k\in\Bbb R$ be a strictly increasing time sequence of impulsive actions with $\lim\limits_{k\to\infty}\tau_k=\infty$.		
	Let $\mathcal U_2\subset\{f\,:\,\Bbb Z_+\to\mathcal U_2\}$ be the space of input signals normed by
	$\|\cdot\|_{\mathcal U_2}$ and taking values in a nonempty subset $U_2$ of some linear normed space. Let
	$g\,:\,X\times U_2\to X$ be a mapping defining impulsive actions and the mapping $\phi$ be defined for all $(t,t_0,x,d_1,d_2)\in\Bbb R\times\Bbb R\times X\times\mathcal U_1\times\mathcal U_2$, $t\ge t_0$.
	
	The following data $\Sigma=(X,\Sigma_c,\mathcal U_2, g,\phi,\mathcal E)$ defines a (forward complete) impulsive system if
	
	$(\Sigma_1)$ for all $(k,x,d_1)\in\Bbb Z_+\times X\times\mathcal U_1$ the system $\Sigma_c$ satisfies
	\begin{equation*}
	\tau_{p(t_0)}-t_0<\epsilon_{t_0,x,d_1},\quad T_k:=\tau_{k+1}-\tau_k<\epsilon_{\tau_k,x,d_1}
	\end{equation*}
	where we denote $p(t_0):=\min\{k\in\Bbb Z_+\,:\tau_k\in\mathcal E_{t_0}\}$ with
	$\mathcal E_{t_0}=[t_0,\infty)\cap\mathcal E$; and
	
	$(\Sigma_2)$  the mapping $\phi$ satisfies
	\begin{equation*}
	\gathered
	\phi(t,t_0,x,d_1,d_2)=\phi_c(t,t_0,x,d_1),\quad\text{ for all }\quad t\in[t_0,\tau_{p(t_0)}],\\
	\phi(t,t_0,x,d_1,d_2)=\phi_c(t,\tau_k,g(\phi(\tau_k,t_0,x,d_1,d_2),d_2(k)),d_1)\\
	\quad \text{ for all }\quad
	t\in(\tau_k,\tau_{k+1}],\quad k\in\Bbb Z_+, k\ge p(t_0).
	\endgathered
	\end{equation*}
\end{definition}
We will denote for short
\begin{equation*}
\gathered
\phi(\tau_k^+,t_0,x,d_1,d_2)=g(\phi(\tau_k,t_0,x,d_1,d_2),d_2(k)),\\
k\ge p(t_0),\; \tau_k\ge t_0
\endgathered
\end{equation*}
The conditions $(\Sigma_c1)$ and $(\Sigma_2)$ imply
\begin{equation*}
\gathered
\lim\limits_{t\to\tau_k+}\phi(t,t_0,x,d_1,d_2)=\phi(\tau_k^+,t_0,x,d_1,d_2),\\
\lim\limits_{t\to\tau_k-}\phi(t,t_0,x,d_1,d_2)=\phi(\tau_k,t_0,x,d_1,d_2);
\endgathered
\end{equation*}
and $(\Sigma_c4)$, $(\Sigma_c5)$, $(\Sigma_2)$ imply that for $t\ge\tau\ge t_0$, $(x,d_1,d_2)\in X\times\mathcal U_1\times\mathcal U_2$
the following holds
\begin{equation}\label{TrPr}
\gathered
\phi(t,t_0,x,d_1,d_2)=\phi(t,\tau,\phi(\tau,t_0,d_1,d_2),d_1,d_2).
\endgathered
\end{equation}
The system $\Sigma_c$ describes the continuous dynamics of the impulsive system $\Sigma$. One can also consider its discrete dynamics separately as a system $\Sigma_d$ defined next
\begin{definition}
	A discrete dynamical system with input $\Sigma_d=(X,g,\phi_d,\mathcal U_2)$ is defined  by a normed state space  $(X,\|\cdot\|_X)$; a space of input signals $\mathcal U_2\subset\{f\,:\,\Bbb Z_+\to U_2\}$ with norm $\|\cdot\|_{\mathcal U_2}$ and values in a nonempty subset
	$U_2$ of a linear normed space; a mapping  $g\,:\,X\times U_2\to  X$; and a mapping
	$\phi_d\,\,:(k,l,x,d_2)\mapsto \phi_d(k,l,x,d_2)$, for $(k,l,x,d_2)\in\Bbb Z_+\times\Bbb Z_+\times X\times\mathcal U_2$, $k\ge l$ such that
	
	$(\Sigma_d1)$  $d(k,k,x,d_2)=x$ and \\ $\phi_d(k+1,l,x,d_2)=g(\phi_d(k,l,x,d_2),d_2(k))$ for all $k\ge l$.
\end{definition}
We assume that $\Sigma$ satisfies the following
\begin{assumption} There exist $\xi,\,\xi_{\tau}\in\mathcal K_{\infty}$, $\tau\in\Bbb R_+$ and $\eta_{\tau},\eta\in\mathcal K_{\infty}$ such that
\begin{equation}\label{1.2}
\gathered
\|\phi_c(t,0,x,d_1)\|\le\xi_{\tau}(\|x\|)+\eta_{\tau}(\|d_1\|_{\mathcal U_1}),\quad t\in[0,\tau],
\endgathered
\end{equation}
where $(x,d_1)\in X\times \mathcal U_1$ and
\begin{equation}\label{1.3}
\gathered
\|g(x,d_2)\|\le\xi(\|x\|)+\eta(\|d_2\|_{\mathcal U_2}),
\endgathered
\end{equation}
where $(x,d_2)\in X\times\mathcal U_2$.
\end{assumption}
Now we define the  main stability property of this paper.
\begin{definition}
For a fixed time sequence $\mathcal E$ of impulsive actions the system $\Sigma$ is called input-to-state stable (ISS) if there
exist functions $\beta_{t_0}\in\mathcal K\mathcal L$, $\gamma_{t_0}\in\mathcal K_{\infty}$, such that
for all $x\in X$ and all $(d_1,d_2)\in \mathcal U_1\times\mathcal U_2$ it holds that
\begin{equation}\label{2.1}
\gathered
\|\phi(t, t_0,x,d_1,d_2)\|_X \le\beta_{t_0}(\|x\|_X,t)+\gamma_{t_0}(d),\quad t\ge t_0
\endgathered
\end{equation}
where $d:=\max\{\|d_1\|_{\mathcal U_1},\|d_2\|_{\mathcal U_2}\}$.
\end{definition}
\begin{definition}
The Lie derivative of a function $V\,:\,X\to\Bbb R$ is defined by
\begin{equation*}
\gathered
\dot V(x,\xi)=\lim_{t\to 0+}\frac{1}{t}(V(\phi_c(t,0,x,\xi))-V(x)),\quad (x,\xi)\in X\times U_1
\endgathered
\end{equation*}
\end{definition}
In many practical particular cases there are simpler expressions for calculation of the Lie derivative possible, see Remark 2.15 in \cite{MiP20}.

The aim of our work is to establish conditions guaranteeing the ISS property for nonlinear impulsive systems. 
Next we define the class of functions that we will use as Lyapunov functions for studying the ISS property.
\begin{definition}
 A continuous function  $V:\,X\to\Bbb R_+$ is called ISS-Lyapunov function if for some
$\alpha_1$, $\alpha_2\in\mathcal K_{\infty}$ it holds that
\begin{equation}\label{2.2}
\gathered
\alpha_1(\|x\|_X)\le V(x)\le\alpha_2(\|x\|_X),\quad x\in X,
\endgathered
\end{equation}
and there exists a function $W\in C(X,\Bbb R)$, $W(0)=0$ such that the Lie derivatives $\dot V(x,\xi)$ and $\dot W(x,\xi)$ exist for all $(x,\xi)\in X\times U_1$,
the sets
\begin{equation}\label{G+-}
\gathered
G^+=\{x\in X\,:\,W(x)\ge 0\},\\
G^-=\{x\in X\,:\,W(x)\le 0\},\\
\endgathered
\end{equation}
are not empty and for some $\chi\in\mathcal K_{\infty}$ and
$\vph_i\in\mathcal P$, $\psi_i\in\mathcal P$, $i=1,2$ so that
\begin{equation}\label{A1}
\gathered
x\in G^{-},\quad \|x\|_X\ge\max\{\chi(\|\xi\|_{U_1}),\chi(\|\mu\|_{U_2})\}\\
\Rightarrow
\begin{cases}
\dot V(x,\xi)\le-\vph_1(V(x)),\\
V(g(x,\mu))\le\psi_1(V(x))
\end{cases}
\endgathered
\end{equation}
\begin{equation}\label{A2}
\gathered
x\in G^{+},\quad \|x\|_X\ge\max\{\chi(\|\xi\|_{U_1}),\chi(\|\mu\|_{U_2})\}\\
\Rightarrow
\begin{cases}
\dot V(x,\xi)\le\vph_2(V(x)),\\
V(g(x,\mu))\le\psi_2(V(x))
\end{cases}
\endgathered
\end{equation}
and
\begin{equation}\label{W-dot}
\gathered
W(x)=0,\;\xi\ne 0,\; \|x\|\ge\chi(\|\xi\|)\Rightarrow \dot W(x,\xi)>0.
\endgathered
\end{equation}
\end{definition}
This definition defers from the known ones, as for example in  \cite{dash-miron13}, due to the auxiliary function $W$ which is of Chetaev type.
\section{Main result}
Our main results establish conditions guaranteeing the ISS property of $\Sigma$.
\begin{theorem}\label{th1}
Let the impulsive system $\Sigma$ satisfy the Assumption 1 and possesses an ISS-Lyapunov function $V$, satisfying \eqref{2.2}--\eqref{A2} such that for some constants $\theta_1$ and $\theta_2$ ($\theta_1\le\theta_2$) and $\delta>0$ for all $a>0$ holds
\begin{equation}\label{A3}
\gathered
\int\limits_a^{\psi_1(a)}\frac{ds}{\vph_1(s)}\le\theta_1-\delta,
\endgathered
\end{equation}
\begin{equation}\label{A4}
\gathered
\int\limits_{\psi_2(a)}^{a}\frac{ds}{\vph_2(s)}\ge\theta_2+\delta.
\endgathered
\end{equation}
Then for any  $\mathcal E$ such that the dwell-time $T_k=\tau_{k+1}-\tau_k$, $k\in\Bbb Z_+$
satisfies $\theta_1\le T_k\le\theta_2$ the system $\Sigma$ is ISS.
\end{theorem}
The proof of this theorem is split into several steps. Without loss of generality we assume $t_0\le\tau_0$, $p(t_0)=0$.

\begin{proposition}\label{prop1}
Let $\Sigma$ satisfy the Assumtion 1 and for some $\beta_{\tau_0^+}\in\mathcal K\mathcal L$, $\gamma_{\tau_0^+}\in\mathcal K_{\infty}$
its solutions satisfy for all $t>\tau_0$ the inequality
\begin{equation}\label{ISS*}
\gathered
\|\phi(t,t_0,\phi_0,d_1,d_2)\|_X\\
\le\beta_{\tau_0^+}(\|\phi(\tau_0^+,t_0,\phi_0,d_1,d_2)\|_X,t)+\gamma_{\tau_0^+}(d),
\endgathered
\end{equation}
Then $\Sigma$ is ISS.
\end{proposition}
{\bf Proof.} Let us fix any initial state $\phi_0$ and disturbance $d_1,d_2$.
From $\Sigma_2$, \eqref{GrPr} and \eqref{1.2} it follows that for all $t\in[t_0,\tau_0]$ the next estimate holds
\begin{equation}\label{O1}
\gathered
\|\phi(t,t_0,\phi_0,d_1,d_2)\|_X=
\|\phi_c(t-t_0,0,\phi_0,\mathcal S_{t_0}d_1)\|_X\\
\le \xi_{\tau_0-t_0}(\|\phi_0\|_X)+\eta_{\tau_0-t_0}(d).
\endgathered
\end{equation}
Also from \eqref{1.3} and \eqref{CF} we have the next estimate
$
\|\phi(\tau_0^+,t_0,\phi_0,d_1,d_2)\|_X\le
 \xi(\|\phi(\tau_0,t_0,\phi_0,d_1,d_2)\|_X)+\eta(d)$
$ \le
\xi(\xi_{\tau_0-t_0}(\|\phi_0\|_X)+\eta_{\tau_0-t_0}(d))+\eta(d)
\\
\le\widehat{\xi}_{\tau_0-t_0}(\|\phi_0\|_X)+\widehat{\eta}_{\tau_0-t_0}(d),$\\
for some fixed $\widehat{\xi}_{\tau_0-t_0}\in\mathcal K_{\infty}$, $\widehat{\eta}_{\tau_0-t_0}\in\mathcal K_{\infty}$.

Finally, from \eqref{ISS*} and \eqref{CF} we have
\begin{equation}\label{O2}
\gathered
\|\phi(t,t_0,\phi_0,d_1,d_2)\|\le \beta_{\tau_0^+}(\|\phi(\tau_0^+,t_0,\phi_0,d_1,d_2)\|_X,t)\\
+\gamma_{\tau_0^+}(d)\le\beta_{\tau_0^+}(\widehat{\xi}_{\tau_0-t_0}(\|\phi_0\|_X)+\widehat{\eta}_{\tau_0-t_0}(d),t)+\gamma_{\tau_0^+}(d)\\
\le\widehat{\beta}_{\tau_0^+}(\|\phi_0\|_X,t)+\wt{\beta}_{\tau_0^+}(\wt{\eta}_{\tau_0-t_0}(d),\tau_0^+)+\gamma_{\tau_0^+}(d)\\
\le\widehat{\beta}_{\tau_0^+}(\|\phi_0\|_X,t)+\widehat{\gamma}_{\tau_0^+}(d),\quad  t>\tau_0
\endgathered
\end{equation}
for some $\widehat{\beta}_{\tau_0^+}$, $\wt{\beta}_{\tau_0^+}\in\mathcal K\mathcal L$, $\widehat{\gamma}_{\tau_0^+}\in\mathcal K_{\infty}$.
Let us define for $s\ge0$ and $t\ge t_0$
\begin{equation*}
\gathered
b_{t_0}(s,t):=
\begin{cases}
\xi_{\tau_0-t_0}(s),\quad t\in[t_0,\tau_0],\\
\xi_{\tau_0-t_0}(s)e^{-t+\tau_0},\quad t>\tau_0;
\end{cases}\\
\beta_{t_0}(s,t):=\max\{b_{t_0}(s,t),\widehat{\beta}_{\tau_0^+}(s,t)\};
\endgathered
\end{equation*}
\begin{equation*}
\gathered
\gamma_{t_0}(s):=\max\{\widehat{\gamma}_{\tau_0^+}(s),\eta_{\tau_0-t_0}(s)\}.
\endgathered
\end{equation*}
From these definitions follows $\beta_{t_0}\in\mathcal K\mathcal L$, $\gamma_{t_0}\in\mathcal K_{\infty}$ and \eqref{O1} with \eqref{O2} imply the ISS property for the impulsive system $\Sigma$, which proves the proposition.

\begin{lemma}\label{lem1}
Let $V$ be an ISS-Lyapunov function of $\Sigma$ satisfying \eqref{2.2}--\eqref{A2} with $\varphi_1,\varphi_2$ satisfying \eqref{A3}-\eqref{A4}.
Fix any $\phi_0\in X$, $(d_1,d_2)\in\mathcal U_1\times\mathcal U_2$ and $r>\chi(d)$.
If $\|\phi(t,t_0,\phi_0,d_1,d_2)\|_X\ge r$ for all $t\in(\tau_p,\tau_m]$ with some $m>p$, then
\begin{equation}\label{inequality-lem1}
\gathered
F(v(\tau_l^+),v(\tau_p^+))\ge \delta(l-p),\quad p\le l\le m
\endgathered
\end{equation}
where $v(t):=V(\phi(t,t_0,\phi_0,d_1,d_2))$ and\\
for any $s>0,\;q>0$
\begin{equation*}
\gathered
F(s,q):=\int\limits_{s}^q\frac{ds}{\widehat{\vph}(s)},\quad \widehat{\vph}(s):=\min\{\vph_1(s),\vph_2(s),s\} ,
\endgathered
\end{equation*}
so that  $F(s,q)\to\infty$ for $s\to 0+$ for any fixed $q>0$.
\end{lemma}

{\bf Proof.} For $l=p$ the inequality \eqref{inequality-lem1} is trivially satisfied. For $l=p+1$ we consider the solution $\phi(t,t_0,\phi_0,d_1,d_2)$ of $\Sigma$ for $t\in(\tau_p,\tau_{p+1}]$ and either of two possible cases 1) $\phi_0\in G^+$ and 2) $\phi_0\in\intt(G^-)$, recalling that $G^+\cup G^-=X$ and $G^+\cap\intt(G^-)=\emptyset$ by definition.

1) For $\phi_0\in G^+$ we will show that $\phi(t,t_0,\phi_0,d_1,d_2)\in G^+$ for $t\in(\tau_p,\tau_{p+1}]$.
Assume, this is not the case, that is there exists
\begin{equation*}
\gathered
\wt t=\sup\{t\in(\tau_p,\tau_{p+1}]\,:\,\phi(t,t_0,\phi_0,d_1,d_2)\in G^+\}\in(\tau_p,\tau_{p+1})
\endgathered
\end{equation*}
such that $\phi(t,t_0,\phi_0,d_1,d_2)\in G^+$, $t\in(\tau_p,\wt t]$ and $\phi(\wt t,t_0,\phi_0,d_1,d_2)\in\partial G^+$.
Let us denote for short $w(t):=W(\phi(t,t_0,\phi_0,d_1,d_2))$, where $t\mapsto W(\phi(t,t_0,\phi_0,d_1,d_2))$ is absolutely continuous.
Hence $w(\wt t)=0$. By assumptions of the lemma we have $\|\phi(\wt t,t_0,\phi_0,d_1,d_2)\|_X\ge r>\chi(d)$, hence due to \eqref{W-dot} and (4.8) from \cite{JMPW21}
follows $\dot w(t)>0$. This means  that for some $\epsilon>0$ for $t\in(\wt t,\wt t+\epsilon)$ it holds that $w(\wt t)\ge 0$, that is $\phi(\wt t,t_0,\phi_0,d_1,d_2)\in G^+$ contradicting the choice of $\wt t$.

This means that $\phi(t,t_0,\phi_0,d_1,d_2)\in G^+$ and \\ $\|\phi(t,t_0,\phi_0,d_1,d_2)\|_X\ge r>\chi(d)$ for $t\in(\tau_p,\tau_{p+1}]$, hence \eqref{A2} imples
\begin{equation*}
\gathered
\dot v(t)\le\vph_2(v(t)),\quad t\in(\tau_p,\tau_{p+1}].
\endgathered
\end{equation*}
We calculate
\begin{equation*}
\gathered
\int\limits_{v(\tau_p^+)}^{v(\tau_{p+1})}\frac{ds}{\vph_2(s)}=
\int\limits_{\tau_p}^{\tau_{p+1}}\frac{dv(s)}{\vph_2(v(s))}\le
\int\limits_{\tau_p}^{\tau_{p+1}}ds=\tau_{p+1}-\tau_p=T_p\le\theta_2.
\endgathered
\end{equation*}
Setting $a=v(\tau_{p+1})$ in \eqref{A4} we obtain
\begin{equation*}
\gathered
\int\limits_{\psi_2(v(\tau_{p+1}))}^{v(\tau_{p+1})}\frac{ds}{\vph_2(s)}\ge\theta_2+\delta\ge
\int\limits_{v(\tau_p^+)}^{v(\tau_{p+1})}\frac{ds}{\vph_2(s)}+\delta.
\endgathered
\end{equation*}
which implies
\begin{equation*}
\gathered
\int\limits_{\psi_2(v(\tau_{p+1}))}^{v(\tau_p^+)}\frac{ds}{\vph_2(s)}\ge\delta.
\endgathered
\end{equation*}
Due to $\|\phi(\tau_{p+1},t_0,\phi_0,d_1,d_2)\|_X\ge r>\chi(d)$ from
\eqref{A4} follows $v(\tau_{p+1}^+)\le\psi_2(v(\tau_{p+1}))$ which implies
\begin{equation*}
\gathered
\int\limits_{v(\tau_{p+1}^+)}^{v(\tau_p^+)}\frac{ds}{\widehat{\vph}(s)}\ge\int\limits_{v(\tau_{p+1}^+)}^{v(\tau_p^+)}\frac{ds}{\vph_2(s)}\ge\delta
\endgathered
\end{equation*}
or equivalently
\begin{equation}\label{H1}
\gathered
F(v(\tau_{p+1}^+),v(\tau_p^+))\ge\delta.
\endgathered
\end{equation}

2) Now let $\phi_0\in\intt G^-$, then either

(i) $\phi(t,t_0,\phi_0,d_1,d_2)\in\intt G^-$ for $t\in(\tau_p,\tau_{p+1}]$ or

(ii) $\phi(\wt t_1,t_0,\phi_0,d_1,d_2)\in G^+$ for some $\wt t_1\in(\tau_p,\tau_{p+1}]$.

In case (i) from $\|\phi(t,t_0,\phi_0,d_1,d_2)\|_X\ge r>\chi(d)$ and \eqref{A1} it follows that
\begin{equation*}
\gathered
\dot v(t)\le-\vph_1(v(t)),\quad t\in(\tau_p,\tau_{p+1}].
\endgathered
\end{equation*}
which means
\begin{equation*}
\gathered
\int\limits_{v(\tau_p^+)}^{v(\tau_{p+1})}\frac{ds}{\vph_1(s)}=\int\limits_{\tau_p}^{\tau_{p+1}}\frac{dv(s)}{\vph_1(v(s))}\le-(\tau_{p+1}-\tau_p),
\endgathered
\end{equation*}
and hence
\begin{equation*}
\gathered
\theta_1\le T_p=\tau_{p+1}-\tau_p\le \int\limits_{v(\tau_{p+1})}^{v(\tau_p^+)}\frac{ds}{\vph_1(s)}.
\endgathered
\end{equation*}
Setting $a=v(\tau_{p+1})$ in \eqref{A3} we obtain
\begin{equation*}
\gathered
\theta_1\ge\int\limits_{v(\tau_{p+1})}^{\psi_1(v(\tau_{p+1}))}\frac{ds}{\vph_1(s)}+\delta.
\endgathered
\end{equation*}
That is
\begin{equation*}
\gathered
\int\limits_{v(\tau_{p+1})}^{v(\tau_p^+)}\frac{ds}{\vph_1(s)}\ge
\int\limits_{v(\tau_{p+1})}^{\psi_1(v(\tau_{p+1})}\frac{ds}{\vph_1(s)}+\delta,
\endgathered
\end{equation*}
and
\begin{equation*}
\gathered
\int\limits_{\psi_1(v(\tau_{p+1}))}^{v(\tau_p^+)}\frac{ds}{\vph_1(s)}\ge
\delta.
\endgathered
\end{equation*}
From $\phi(\tau_{p+1},t_0,\phi_0,d_1,d_2)\in G^-$ and\\ $\|\phi(\tau_{p+1},t_0,\phi_0,d_1,d_2)\|_X\ge r>\chi(d)$,
follows\\
$v(\tau_{p+1}^+)\le\psi_1(v(\tau_{p+1}))$, and hence
\begin{equation*}
\gathered
\int\limits_{v(\tau_{p+1}^+)}^{v(\tau_p^+)}\frac{ds}{\widehat{\vph}(s)}\ge\int\limits_{v(\tau_{p+1}^+)}^{v(\tau_p^+)}\frac{ds}{\vph_1(s)}\ge \int\limits_{\psi_1(v(\tau_{p+1}))}^{v(\tau_p^+)}\frac{ds}{\vph_1(s)}\ge
\delta
\endgathered
\end{equation*}
or in other words
\begin{equation}\label{H2}
\gathered
F(v(\tau_{p+1}^+),v(\tau_p^+))
\ge\delta.
\endgathered
\end{equation}
In case (ii) we define
\begin{equation*}
\gathered
\widehat t=\inf\{t\in(\tau_p,\tau_{p+1}]\,:\,\phi(t,\tau_p^+,\phi_0,d_1,d_2)\in\intt G^+\},
\endgathered
\end{equation*}
so that  $\phi(\widehat t,t_0,\phi_0,d_1,d_2)\in\partial G^+\subset G^+$  and $\widehat t>\tau_p$. From the properties of  $W$ it follows that
$\phi(t,t_0,\phi_0,d_1,d_2)\in\intt G^-$ for $t\in(\tau_{p},\widehat t)$ and
$\phi(t,t_0,\phi_0,d_1,d_2)\in G^+$ for $t\in[\widehat t,\tau_{p+1}]$ (similarly to the case 1) above).

From \eqref{A1} follows
\begin{equation}\label{*}
\gathered
v(\widehat t)\le v(\tau_{p}^+).
\endgathered
\end{equation}
Since
$\|\phi(t,t_0,\phi_0,d_1,d_2)\|_X\ge r>\chi(d)$ and\\
$\phi(t,t_0,\phi_0,d_1,d_2)\in G^+$ for
$t\in[\widehat t,\tau_{p+1}]$, then from \eqref{A2} follows
\begin{equation*}
\gathered
\dot v(t)\le\vph_2(v(t)),\quad t\in[\widehat t,\tau_{p+1}].
\endgathered
\end{equation*}
This allows to calculate
\begin{equation*}
\gathered
\int\limits_{v(\widehat t)}^{v(\tau_{p+1})}\frac{ds}{\vph_2(s)}=
\int\limits_{\widehat t}^{\tau_{p+1}}\frac{dv(s)}{\vph_2(v(s))}\le\tau_{p+1}-\widehat t\le\tau_{p+1}-\tau_p= T_p\le\theta_2,
\endgathered
\end{equation*}
so that \eqref{*} implies
\begin{equation*}
\gathered
\int\limits_{v(\tau_p^+)}^{v(\tau_{p+1})}\frac{ds}{\vph_2(s)}\le\int\limits_{v(\widehat t)}^{v(\tau_{p+1})}\frac{ds}{\vph_2(s)}\le\theta_2.
\endgathered
\end{equation*}
Now we set $a=v(\tau_p)$ into \eqref{A4} and obtain
\begin{equation*}
\gathered
\int\limits_{v(\tau_p^+)}^{v(\tau_{p+1})}\frac{ds}{\vph_2(s)}\le\theta_2\le\int\limits_{\psi_2(v(\tau_{p+1}))}^{v(\tau_{p+1})}\frac{ds}{\vph_2(s)}-\delta
\endgathered
\end{equation*}
or
\begin{equation*}
\gathered
\int\limits_{v(\tau_p^+)}^{\psi_2(v(\tau_{p+1}))}\frac{ds}{\vph_2(s)}\le-\delta.
\endgathered
\end{equation*}
Having $\|\phi(\tau_{p+1},t_0,\phi_0,d_1,d_2)\|_X\ge r>\chi(d)$,\\ $\phi(\tau_{p+1},t_0,\phi_0,d_1,d_2)\in G^+$ and \eqref{A4} we conclude that $v(\tau_{p+1}^+)\le\psi_2(v(\tau_{p+1}))$. Hence
\begin{equation*}
\gathered
\int\limits_{v(\tau_p^+)}^{v(\tau_{p+1}^+)}\frac{ds}{\widehat{\vph}(s)}\le\int\limits_{v(\tau_p^+)}^{v(\tau_{p+1}^+)}\frac{ds}{\vph_2(s)}\le
\int\limits_{v(\tau_p^+)}^{\psi_2(v(\tau_{p+1}))}\frac{ds}{\vph_2(s)}\le-\delta.
\endgathered
\end{equation*}
In other words
\begin{equation}\label{H3}
\gathered
F(v(\tau_{p+1}^+),v(\tau_p^+))\ge\delta.
\endgathered
\end{equation}
From \eqref{H1}---\eqref{H3} we conclude that \eqref{inequality-lem1} is true for $l=p+1$. Considering the solution between the next two consequent jumps we obtain
\begin{equation}\label{subintervals}
\gathered
F(v(\tau_{p+1}^+),v(\tau_l^+))\ge\delta,\quad l=p,p+1,\dots,m-1.
\endgathered
\end{equation}
From the definition of $F$ for any $s\ge z\ge q>0$ we have $F(s,z)+F(z,q)=F(s,q)$. Since the interval $(\tau_p,\tau_l)$ is split into $l-p$ subintervals by the time instants of the impulsive actions we finally obtain from \eqref{subintervals} that \eqref{inequality-lem1} is proved.

\begin{remark}\label{rem2}
 Let $F^{-1}(s,\cdot)$ be the inverse function to $F(\cdot,s)$, $s\in\Bbb R_+$. From $F(\tau,s)\to\infty$ for $\tau\to 0+$ follows $F^{-1}(s,\tau)\to 0$ for $\tau\to+\infty$. Also note that $F^{-1}(s,\cdot)$ is strictly decreasing whereas $F^{-1}(\cdot,s)$ is strictly increasing for $s>0$.
\end{remark}

\begin{lemma}
Under the conditions of Theorem \ref{th1} let $r$ be such that $r>\chi(d)$, then there exists $\tau_{k_0}\ge\tau_0$ such that $\phi(\tau_{k_0}^+,t_0,\phi_0,d_1,d_2)\in B_r(0)$.
\end{lemma}

{\bf Proof.} Assume by contradiction that for all $k\in\Bbb Z_+$ we have $\|\phi(\tau_k,t_0,\phi_0,d_1,d_2)\|_X\ge r$. Lemma  \ref{lem1} implies that the sequence
$\{v(\tau_k^+)\}_{k=0}^{\infty}$ is decreasing. Being bounded from below it has a limit $v^*\ge 0$.
From \eqref{inequality-lem1} we have
\begin{equation*}
\gathered
F(v(\tau_{m}^+),v(\tau_0^+))\ge\delta m.
\endgathered
\end{equation*}
If $v^*\ne 0$ this inequality leads to a contradiction  letting $m\to\infty$.
Hence, $v^*=0$.  Due to \eqref{2.2} we have
\begin{equation*}
\gathered
0\le \alpha_1(r)\le\alpha_1(\|\phi(\tau_m^+,t_0,\phi_0,d_1,d_2)\|_X)\le\\ V(\phi(\tau_m^+,t_0,\phi_0,d_1,d_2))=v(\tau_m^+)
\endgathered
\end{equation*}
and taking the limit for  $m\to\infty$ we arrive to $\alpha_1(r)=0$ which implies $r=0$ contradicting $r>\chi(d)$. This finishes the proof of the lemma.
\begin{lemma}\label{lem3}
Under the conditions of Theorem \ref{th1} the solution $\phi$ of $\Sigma$ satisfy
$\phi(\tau_{k_0}^+,t_0,\phi_0,d_1,d_2)\in B_r(0)$ for some $\tau_{k_0}\in{\mathcal E}$ and $r>\chi(d)$, then there exists $R\in\mathcal K^2_{\infty}$ such that
$$\|\phi(t,t_0,\phi_0,d_1,d_2)\|\le R(r,d),\quad  t>\tau_{k_0}.$$
\end{lemma}
{\bf Proof.} We define for $s\ge0,\;q\ge0$
\begin{equation*}
\gathered
R(s,q):=\max\{R_1(s,q),R_4(s,q),R_6(s,q),s\},
\endgathered
\end{equation*}
where we use the following combinations of functions from \eqref{1.2}-\eqref{1.3}
\begin{equation*}
\gathered
R_1(s,q)=\xi_{\theta_2}(s)+\eta_{\theta_2}(q),\quad R_2(s,q)=\xi(s)+\eta(q),\\
R_3(s,q)=\max\{\xi(R_1(s,q))+\eta(q),R_2(s,q)\},\\ R_4(s,q)=\xi_{\theta_2}(R_3(s,q))+\eta_{\theta_2}(q),\\
R_5(s,q)=(\alpha_1^{-1}\circ\alpha_2)(R_3(s,q)),\\
R_6(s,q)=\xi_{\theta_2}(R_5(s,q))+\eta_{\theta_2}(q).
\endgathered
\end{equation*}
Let $\widehat t_1>\tau_{k_0}$ be such that
$\|\phi(\widehat t_1,t_0,\phi_0,d_1,d_2)\|_X\le r$
and
$\|\phi(\widehat t_1^+,t_0,\phi_0,d_1,d_2)\|_X\ge r$
(if such  $\widehat t_1$ does not exist, then the result is proved).
Define
$$\gathered\widehat t_2:=\sup\{t>\widehat t_1:\,\|\phi(s,t_0,\phi_0,d_1,d_2)\|_X\ge r\\
\text{for}\quad s\in[\widehat t_1,t]\}\in[\widehat t_1,\infty],\endgathered$$
so that
\begin{equation*}
\gathered
\|\phi(t,t_0,\phi_0,d_1,d_2)\|_X\ge r, \quad t\in[\widehat t_1,\widehat t_2].
\endgathered
\end{equation*}
It is enough to show that $\|\phi(t,t_0,\phi_0,d_1,d_2)\|_X\le R(r,d)$ for $t\in[\widehat t_1,\widehat t_2]$.

If $\mathcal E_{[\widehat t_1,\widehat t_2]}:=[\widehat t_1,\widehat t_2]\cap\mathcal E=\emptyset$, then by the properties \eqref{TrPr}, $\Sigma_2$,
\eqref{GrPr} for all $t\in[\widehat t_1,\widehat t_2]$ we conclude
\begin{equation*}
\gathered
\|\phi(t,t_0,\phi_0,d_1,d_2)\|_X\\
=\|\phi(t,\widehat t_1,\phi(\widehat t_1,t_0,\phi_0,d_1,d_2),d_1,d_2)\|_X\\
=\|\phi_c(t,\widehat t_1,\phi(\widehat t_1,t_0,\phi_0,d_1,d_2),d_1)\|_X\\=
\|\phi_c(t-\widehat t_1,0,\phi(\widehat t_1,t_0,\phi_0,d_1,d_2),\mathcal S_{\widehat t_1}d_1)\|_X\\
\le\xi_{\theta_2}(r)+\eta_{\theta_2}(d)=R_1(r,d).
\endgathered
\end{equation*}
If otherwise $\mathcal E_{[\widehat t_1,\widehat t_2]}\ne\emptyset$ we denote its minimal element by $\tau_p\ge\widehat t_1$ and consider two possible cases (i) $\widehat t_1=\tau_p$ and (ii) $\tau_p>\widehat t_1$ separately.

In case of (i) from \eqref{1.3} follows
\begin{equation*}
\gathered
\|\phi(\widehat t_1^+,t_0,\phi_0,d_1,d_2)\|_X\le \xi(\|\phi(\widehat t_1,t_0,\phi_0,d_1,d_2)\|_X)\\+\eta(\|d_2\|_{\mathcal U_2})
\le
\xi(r)+\eta(d)=R_2(r,d).
\endgathered
\end{equation*}
In case of (ii) by means of \eqref{1.2} and \eqref{1.3} we obtain
\begin{equation*}
\gathered
\|\phi(\tau_p^+,t_0,\phi_0,d_1,d_2)\|_X\le\xi(\|\phi(\tau_p,t_0,\phi_0,d_1,d_2)\|_X)+\eta(\|d_2\|_{\mathcal U_2})
\endgathered
\end{equation*}
and with help of \eqref{TrPr}, $\Sigma_2$, \eqref{GrPr} we get
\begin{equation*}
\gathered
\|\phi(\tau_p,t_0,\phi_0,d_1,d_2)\|_X\\
=\|\phi(\tau_p,\widehat t_1,\phi(\widehat t_1,t_0,\phi_0,d_1,d_2),d_1,d_2)\|_X\\
=\|\phi_c(\tau_p,\widehat t_1,\phi(\widehat t_1,t_0,\phi_0,d_1,d_2),d_1)\|_X\\=
\|\phi_c(\tau_p-\widehat t_1,0,\phi(\widehat t_1,t_0,\phi_0,d_1,d_2),\mathcal S_{\widehat t_1}d_1)\|_X\\
\le\xi_{\theta_2}(\|\phi(\widehat t_1,t_0,\phi_0,d_1,d_2)\|_X)+\eta_{\theta_2}(\|d_1\|_{\mathcal U_2})\\
\le\xi_{\theta_2}(r)+\eta_{\theta_2}(d)=R_1(r,d).
\endgathered
\end{equation*}
Hence,
\begin{equation*}
\gathered
\|\phi(\tau_p^+,t_0,\phi_0,d_1,d_2)\|_X\le \xi(R_1(r,d))+\eta(d).
\endgathered
\end{equation*}
In bothe cases (i) and (ii) we see that \\
$\|\phi(\tau_p^+,t_0,\phi_0,d_1,d_2)\|_X\le R_3(r,d)$ holds.

If $\sharp\mathcal E_{[\widehat t_1,\widehat t_2]}=1$, then by means of $\Sigma_2$, \eqref{GrPr}, \eqref{1.2} and \eqref{1.3} we obtain for
$t\in(\tau_p,\widehat t_2]$
\begin{equation*}
\gathered
\|\phi(t,t_0,\phi_0,d_1,d_2)\|_X\\
=\|\phi_c(t-\tau_p,0,\phi(\tau_p^+,t_0,\phi_0,d_1,d_2),\mathcal S_{\tau_p}d_1)\|_X\\
\le\xi_{\theta_2}(\|\phi(\tau_p^+,t_0,\phi_0,d_1,d_2)\|_X)+\eta_{\theta_2}(\|d_1\|_{\mathcal U_1})\\
\le\xi_{\theta_2}(R_3(r,d))+\eta_{\theta_2}(d)=R_4(r,d).
\endgathered
\end{equation*}
If $\sharp\mathcal E_{[\widehat t_1,\widehat t_2]}\ge 2$, then by Lemma \ref{lem1} and observing that $F(s,q)>0\;\Leftrightarrow\;s<q$ we obtain for $l\ge p$ such that $\tau_l\in\mathcal E_{[\widehat t_1,\widehat t_2]}$
\begin{equation*}
\gathered
\alpha_1(\|\phi(\tau_l^+,t_0,\phi_0,d_1,d_2)\|_X)\le V(\phi(\tau_l^+,t_0,\phi_0,d_1,d_2))\\
=v(\tau_l^+)\le v(\tau_p^+)=V(\phi(\tau_p^+,t_0,\phi_0,d_1,d_2))\\
\le \alpha_2(\|\phi(\tau_p^+,t_0,\phi_0,d_1,d_2)\|_X)\le\alpha_2(R_3(r,d)).
\endgathered
\end{equation*}
This implies $\|\phi(\tau_l^+,t_0,\phi_0,d_1,d_2)\|_X\le(\alpha_1^{-1}\circ\alpha_2)(R_3(r,d))\\
=R_5(r,d)$.
Hence for all $t\in(\tau_l,\tau_{l+1}]$, from \eqref{1.2} and properties $\Sigma_2$, \eqref{GrPr}, \eqref{TrPr} we obtain
\begin{equation*}
\gathered
\|\phi(t,t_0,\phi_0,d_1,d_2)\|_X\\
=\|\phi_c(t,\tau_l,\phi(\tau_l^+,t_0,\phi_0,d_1,d_2),d_1)\|_X\\=
\|\phi_c(t-\tau_l,0,\phi(\tau_l^+,t_0,\phi_0,d_1,d_2),\mathcal S_{\tau_l}d_1)\|_X\\
\le\xi_{\theta_2}(\|\phi(\tau_l^+,t_0,\phi_0,d_1,d_2)\|_X)+\eta_{\theta_2}(\|d_1\|_{\mathcal U_1})\\
\le\xi_{\theta_2}(R_5(r,d))+\eta_{\theta_2}(d)=R_6(r,d).
\endgathered
\end{equation*}
This implies the estimate
\begin{equation*}
\gathered
\|\phi(t,t_0,\phi_0,d_1,d_2)\|_X\le R(r,d)\quad\text{for all}\quad t\in[\tau_p,\tau_{m+1}].
\endgathered
\end{equation*}
If $t\in[\widehat t_1,\tau_p]$, then from $\Sigma_2$, \eqref{GrPr}, \eqref{TrPr} and \eqref{1.2} we get
\begin{equation*}
\gathered
\|\phi(t,t_0,\phi_0,d_1,d_2)\|_X\\
=\|\phi(t,\widehat t_1,\phi(\widehat t_1,t_0,\phi_0,d_1,d_2),d_1,d_2)\|_X\\
=\|\phi_c(t,\widehat t_1,\phi(\widehat t_1,t_0,\phi_0,d_1,d_2),d_1)\|_X\\
=\|\phi_c(t-\widehat t_1,0,\phi(\widehat t_1,t_0,\phi_0,d_1,d_2),\mathcal S_{\widehat t_1}d_1)\|_X\\
\le\xi_{\theta_2}(\|\phi(\widehat t_1,t_0,\phi_0,d_1,d_2)\|_X)+\eta_{\theta_2}(d)\\
=\xi_{\theta_2}(r)+\eta_{\theta_2}(d)=R_1(r,d)\le R(r,d).
\endgathered
\end{equation*}
Since $[\widehat t_1,\tau_{m+1}]\supseteq[\widehat t_1,\widehat t_2]$ holds, the lemma is proved.

{\bf Proof of Theorem \ref{th1}} Take $r=(1+\varepsilon)\chi(d)$ for some $\varepsilon>0$ and denote by $k_0$ the smallest integer for which $\phi(\tau_{k_0}^+,t_0,\phi_0,d_1,d_2)\in B_r(0)$ holds, that is $\|\phi(\tau_{k}^+,t_0,\phi_0,d_1,d_2)\|_X\ge r$ for all $0\le k\le k_0-1$.
From Lemma \ref{lem1} follows $F(v(\tau_k^+),v(\tau_0^+))\ge\delta k$, hence (see Remark \ref{rem2})  we have
\begin{equation*}
\gathered
v(\tau_k^+)\le F^{-1}(v(\tau_0^+),k\delta).
\endgathered
\end{equation*}
From \eqref{1.2} and properties $\Sigma_2$, \eqref{TrPr} and \eqref{GrPr} it follows that for $t\in(\tau_k,\tau_{k+1}]$, $k=0,\dots,k_0-1$ the next inequality holds
\begin{equation*}
\gathered
\|\phi(t,t_0,\phi_0,d_1,d_2)\|_X=\|\phi_c(t,\tau_k,\phi(\tau_k^+,t_0,\phi_0,d_1,d_2),d_1,d_2)\|_X\\
=\|\phi_c(t-\tau_k,0,\phi(\tau_k^+,t_0,\phi_0,d_1,d_2),\mathcal S_{\tau_k}d_1)\|_X\\
\le
\xi_{\theta_2}(\|\phi(\tau_k^+,t_0,\phi_0,d_1,d_2)\|_X)+\eta_{\theta_2}(d).
\endgathered
\end{equation*}
Now \eqref{2.2} implies
\begin{equation*}
\gathered
\alpha_1(\|\phi(\tau_k^+,t_0,\phi_0,d_1,d_2)\|_X)\le V(\phi(\tau_k^+,t_0,\phi_0,d_1,d_2))\\
=v(\tau_k^+)
\le F^{-1}(v(\tau_0^+),\delta k),\\
v(\tau_0^+)=V(\phi(\tau^+_0,t_0,\phi_0,d_1,d_2))\le\alpha_2(\|\phi(\tau^+_0,t_0,\phi_0,d_1,d_2)\|_X).
\endgathered
\end{equation*}
Defining $\beta_k(s):=F^{-1}(\alpha_2(s),\delta k),\;s>0$ and $\beta_k(0):=0$ we can write
\begin{equation*}
\gathered
\alpha_1(\|\phi(\tau_k^+,t_0,\phi_0,d_1,d_2)\|_X)\\
\le F^{-1}(\alpha_2(\|\phi(\tau^+_0,t_0,\phi_0,d_1,d_2)\|_X),\delta k)\\
=\beta_k(\|\phi(\tau^+_0,t_0,\phi_0,d_1,d_2)\|_X).
\endgathered
\end{equation*}
Which means
\begin{equation*}
 \|\phi(\tau_k^+,t_0,\phi_0,d_1,d_2)\|_X\le(\alpha_1^{-1}\circ\beta_k)(\|\phi(\tau^+_0,t_0,\phi_0,d_1,d_2)\|_X).
 \end{equation*}
 By the properties of $F^{-1}$ (see Remark \ref{rem2}) it follows that $\beta_k\in\mathcal K$, $\beta_{k+1}(s)<\beta_k(s)$ and $\lim_{k\to\infty}\beta_k(s)=0$ for $s>0$.

For $t\in(\tau_k,\tau_{k+1}]$ the following estimate holds
\begin{equation*}
\gathered
\|\phi(t,t_0,\phi_0,d_1,d_2)\|_X\\ \le
(\xi_{\theta_2}\circ\alpha_1^{-1}\circ\beta_k)(\|\phi(\tau^+_0,t_0,\phi_0,d_1,d_2)\|_X)
+\eta_{\theta_2}(d).
\endgathered
\end{equation*}
Defining $\widehat\beta_k(s):=(\xi_{\theta_2}\circ\alpha_1^{-1}\circ\beta_k)(s)$ for  $k\in\Bbb Z_+$ and $s\ge0$ we see that
$\widehat\beta_k\in\mathcal K$, $\widehat\beta_{k+1}(s)<\widehat\beta_k(s)$ and  $\lim_{k\to\infty}\widehat{\beta}_k(s)=0$ for any $s>0$.
Further, we define
\begin{equation*}
\gathered
\beta_{\tau_0^+}(s,t)=\widehat\beta_k(s)+\frac{t-\tau_k}{\tau_{k+1}-\tau_k}(\widehat\beta_{k+1}(s)-\widehat\beta_k(s)),\\
t\in(\tau_k,\tau_{k+1}],\;s\ge0
\endgathered
\end{equation*}
so that $\beta_{\tau_0^+}\in\mathcal K\mathcal L$ and for $t\in(\tau_0,\tau_{k_0}]$ the following estimate holds
\begin{equation}\label{ISS1}
\gathered
\|\phi(t,t_0,\phi_0,d_1,d_2)\|_X\le
\beta_{\tau_0^+}(\|\phi(\tau^+_0,t_0,\phi_0,d_1,d_2)\|_X,t)\\
+\eta_{\theta_2}(d).
\endgathered
\end{equation}
From Lemma \ref{lem3} follows that for any $t>\tau_{k_0}$ we have
\begin{equation}\label{ISS2}
\gathered
\|\phi(t,t_0,\phi_0,d_1,d_2)\|_X\le \sigma(d),
\endgathered
\end{equation}
where we denote $\sigma(s):=R((1+\varepsilon)s,s)$, $\varepsilon>0$. Since $\sigma\in\mathcal K$ then collecting the estimates \eqref{ISS1} and \eqref{ISS2} we obtain
\begin{equation*}
\gathered
\|\phi(t,t_0,\phi_0,d_1,d_2)\|_X\le
\beta_{\tau_0^+}(\|\phi(\tau^+_0,t_0,\phi_0,d_1,d_2)\|_X,t)\\
+\gamma_{\tau_0^+}(d),\quad t>\tau_0,
\endgathered
\end{equation*}
where $\gamma_{\tau_0^+}(s):=\sigma(s)+\eta_{\theta_2}(s)$.
The statement of the theorem follows then from the Proposition \ref{prop1}.
\begin{remark}
In the definition of an ISS-Lyapunov function $V$ we have assumed that the sets $G^+$ and $G^-$ are not empty. Observe however that in case $G^+=\emptyset$, $G^-=X$ Theorem \ref{th1} reduces to Theorem 1 from \cite{dash-miron13} and if $G^-=\emptyset$, $G^+=X$ then Theorem \ref{th1} reduces to Theorem 3 from  \cite{dash-miron13}.
\end{remark}
\section{Nonlinear interconnection of an ODE and a PDE}
Let $\widehat X=L^2[0,l]\times\Bbb R$, $X=H_0^1(0,l)\times\Bbb R$,
$U_1=U_2=H_0^1(0,l)\times\Bbb R$ and the spaces of input signals be
\begin{equation*}
\gathered
\mathcal U_1=L^{\infty}(U_1)\cap(H_{\loc}(\Bbb R_+,H_0^1(0,l))\times H_{\loc}(\Bbb R_+,\Bbb R)),\\
\mathcal U_2=l^{\infty}(\Bbb Z_+,U_2).
\endgathered
\end{equation*}
Consider the following nonlinear impulsive system
\begin{equation}\label{Exampl}
\gathered
\partial_t x(z,t)=a^2\partial_{zz} x(z,t)+\Phi(x(z,t))\\
+B(z)y(t)+d_{11}(z,t),\quad t\ne\tau_k,\\
\dot y(t)=c^2 y(t)+\int\limits_0^lD(z)x(z,t)\,dz+d_{12}(t),\quad t\ne\tau_k,\\
x(z,t^+)=\alpha(z)x(z,t)+\beta(z)y(t)+d_{21}(z,k),\quad t=\tau_k,\\
y(t^+)=\int\limits_{0}^l\gamma(z)x(z,t)\,dz+\delta y(t)+d_{22}(k),\quad t=\tau_k,
\endgathered
\end{equation}
with initial and boundary conditions
\begin{equation}\label{ExamplBC}
\gathered
x(z,0)=x_0(z)\in H_0^1(0,l),\quad z\in[0,l],\\
\quad y(0)=y_0\in\Bbb R,\\
 x(0,t)=x(l,t)=0,\quad t\in\Bbb R_+,\\
 x_0(0)=x_0(l)=0.
\endgathered
\end{equation}
Here $a$, $c$, $l$ are given positive constants, $\alpha\in C^2[0,l]$, $\,\beta\in H^1_0(0,l)$, $\gamma\in L^2[0,l]$, $D\in L^2[0,l]$, $B\in H_0^1(0,l)$ are given functions and
$d_1(t)=(d_{11}(\cdot,t),d_{12}(t))\in U_1$, $d_2(k)=(d_{21}(\cdot,k),d_{22}(k))\in U_2$ are unknown disturbances.

Assume that $\,\Phi\,:\,\Bbb R\to\Bbb R$ satisfies:

(i) $\Phi\in C^1(\Bbb R)$ with locally Lipshitz $\Phi^{\prime}$ that is for any $s_0\in\Bbb R$ and $\varrho_0>0$ there exists $L=L(s_0,\varrho_0)>0$ such that for all $s\in\Bbb R$ with $|s-s_0|\le\varrho_0$ it holds that $|\Phi^{\prime}(s)-\Phi^{\prime}(s_0)|\le L|s-s_0|$.

(ii) It holds that $\Phi^{\prime}(s)\le 0$, $s\Phi(s)\le 0$ for all $s\in\Bbb R$.

The problem  \eqref{Exampl}--\eqref{ExamplBC} can be written in the following form
\begin{equation}\label{Exampl**}
\gathered
\frac{d}{dt}
\begin{pmatrix}
x(\cdot,t)\\
y(t)
\end{pmatrix}+\mathcal A\begin{pmatrix}
x(\cdot,t)\\
y(t)
\end{pmatrix}=f(t,x,y,d_1),\quad t\ne\tau_k\\
\begin{pmatrix}
x(\cdot,t^+)\\
y(t^+)
\end{pmatrix}=\mathcal B
\begin{pmatrix}
x(\cdot,t)\\
y(t)
\end{pmatrix}+d_2(k),\quad t=\tau_k,
\endgathered
\end{equation}
where $\mathcal A$ is the linear operator on $\widehat X$ defined by
\begin{equation*}
\gathered
\mathcal A\begin{pmatrix}
x\\
y
\end{pmatrix}=\begin{pmatrix}
-a^2\partial_{zz}x(z)\\
-c^2y
\end{pmatrix}
\endgathered
\end{equation*}
with domain $\mathcal D(\mathcal A)=(H_0^1(0,l)\cap H^2(0,l))\times\Bbb R$, $\mathcal B$ is the linear operator defined by
\begin{equation*}
\gathered
 \mathcal B\begin{pmatrix}
x\\
y
\end{pmatrix}=\begin{pmatrix}
\alpha(z)x(z)+\beta(z)y\\
\int\limits_{0}^l\gamma(z)x(z)\,dz+\delta y
\end{pmatrix}
\endgathered
\end{equation*}
and $d_2(k)=(d_{21}(\cdot,k),d_{22}(k))$.

The operator $\mathcal A$ being a direct product of two sectorial operators is sectorial and hence it generates on
$\widehat X$ an analytic semi-group \cite{henri}. The mapping $f\,:\,\Bbb R\times H_0^1(0,l)\times \Bbb R \to L^2[0,l]\times \Bbb R$ in \eqref{Exampl**} is defined by
\begin{equation*}
\gathered
f(t,x,y,d_1)=
\begin{pmatrix}
\Phi(x(z))+B(z)y+d_{11}(z,t)\\
\int\limits_0^lD(z)x(z)\,dz+d_{12}(t)
\end{pmatrix}
\endgathered
\end{equation*}
We consider classical solutions of \eqref{Exampl}---\eqref{ExamplBC} defined as in Definition 3.3.1 of \cite{henri}.
\begin{remark}
For any initial state $(x_0,y_0)\in H_0^1(0,l)\times\Bbb R$ and input $(d_1,d_2)\in\mathcal U_1\times\mathcal U_2$ there exists a unique solution to the propblem \eqref{Exampl}--\eqref{ExamplBC}.
This follows from the fact the the corresponding problem without impulsive actions
\begin{equation}\label{ExamplCD}
\gathered
\frac{d}{dt}
\begin{pmatrix}
x(\cdot,t)\\
y(t)
\end{pmatrix}+\mathcal A\begin{pmatrix}
x(\cdot,t)\\
y(t)
\end{pmatrix}=f(t,x,y),
\endgathered
\end{equation}
for each $(x_0,y_0)\in H_0^1(0,l)\times \Bbb R$ possesses a unique solution defined for $[t_0,t_0+\epsilon_{t_0,x_0,y_0,d_1}]$,
$\epsilon_{t_0,x_0,y_0,d_1}>\theta_2$,
and that the mapping $g(x,y,\mu)=\mathcal B(x,y)+\mu$ keeps the space $X$ invariant  for any $\mu\in U_2$.
\end{remark}

The well-posedness of \eqref{ExamplCD} can be established by means of Theorems 3.3.3 and 3.3.4 from \cite{henri}.
In order to check the conditions of Theorem 3.3.3 (setting $\alpha=0.5$ there) we need to show that the mapping $f$ is locally H\"older wrt $t$ and locally Lipschitz
wrt $(x,y)\in H_0^1(0,l)\times\Bbb R$, that is for any $(t_0,x_0,y_0)\in\Bbb R\times H_0^1(0,l)\times\Bbb R$
there exists
\begin{equation*}
\gathered
O_{\varrho}(t_0,x_0,y_0):=\{(t,x,y)\in\Bbb R_+\times H_0^1(0,l)\times\Bbb R\,\,|\\
\,\,|t-t_0|<\varrho,\, \|x-x_0\|_{H_0^1(0,l)}<\varrho,\, |y-y_0|<\varrho\}.
\endgathered
\end{equation*}
such that for any point $(t,x,y)\in O_{\varrho}(t_0,x_0,y_0)$ it holds that
\begin{equation*}
\gathered
\|f(t,x,y)-f(t_0,x_0,y_0)\|_{\widehat X}\le L(|t-t_0|^{\nu_1}\\
+\|x-x_0\|_{H_0^1(0,l)}+|y-y_0|)
\endgathered
\end{equation*}
for some constants $L>0$, $\nu_1\in(0,1]$.

Indeed,
\begin{equation}\label{loclip}
\gathered
\|f(t,x,y)-f(t_0,x_0,y_0)\|_{\widehat X}\\
\le\frac{l}{\pi} \|\Phi\circ x-\Phi\circ x_0+B(y-y_0)+d_{11}(\cdot,t)-d_{11}(\cdot,t_0)\|_{H_0^1(0,l)}\\
+\Big|\int\limits_0^l D(z)(x(z)-x_0(z))\,dz+d_{12}(t)-d_{12}(t_0)\Big|\\
\le\frac{l}{\pi}(\|\Phi\circ x-\Phi\circ x_0\|_{H_0^1(0,l)}+\|B\|_{H^1_0(0,l)}|y-y_0|)\\
+\|D\|_{L^2[0,l]}\|x-x_0\|_{L^2[0,l]}+|d_{12}(t)-d_{12}(t_0)|\\
+\frac{l}{\pi}\|d_{11}(\cdot,t)-d_{11}(\cdot,t_0)\|_{H^1_0(0,l)}.
\endgathered
\end{equation}
Consider the first summand separately
\begin{equation*}
\gathered
\|\Phi\circ x-\Phi\circ x_0\|_{H^1_0(0,l)}=\|\partial_z(\Phi\circ x-\Phi\circ x_0)\|_{L^2[0,l]}\\
=\|(\Phi^{\prime}\circ x)\, \partial_zx-(\Phi^{\prime}\circ x_0)\, \partial_zx_0\|_{L^2[0,l]}\\
\le \|(\Phi^{\prime}\circ x-\Phi^{\prime}\circ x_0)\partial_z x\|_{L^2[0,l]}\\
+\|(\Phi^{\prime}\circ x_0)(\partial_z x-\partial_z x_0)\|_{L^2[0,l]}
\endgathered
\end{equation*}
By the Sobolev embedding theorem we have  $x\in C[0,l]$ and satisfies
$\|x-x_0\|_{C[0,l]}\le C_1\|x-x_0\|_{H^1_0(0,l)}\le C_1\varrho$, for some $C_1>0$.
Hance, using the condition (i), we obtain
\begin{equation*}
\gathered
|\Phi^{\prime}(x(z))-\Phi^{\prime}(x_0(z))|\le L_1|x(z)-x_0(z)|\\
\le L_1\|x-x_0\|_{C[0,l]}\le L_1C_1\|x-x_0\|_{H^1_0(0,l)},
\endgathered
\end{equation*}
where $L_1$ is a positive constant, which can depend on  $x_0$ and $\rho$. Hence
\begin{equation*}
\gathered
\|\Phi^{\prime}\circ x-\Phi^{\prime}\circ x_0\|_{C[0,l]}\le C_2\|x-x_0\|_{H^1_0(0,l)}\le C_2\varrho,
\endgathered
\end{equation*}
where $C_2=L_1C_1$, and
\begin{equation*}
\gathered
\|\Phi\circ x-\Phi\circ x_0\|_{H^1_0(0,l)}\\
\le \|\Phi^{\prime}\circ x-\Phi^{\prime}\circ x_0\|_{C[0,l]}\|\partial_z x\|_{L^2[0,l]}\\
+\|\Phi^{\prime}\circ x_0\|_{C[0,l]}\|\partial_z x-\partial_z x_0\|_{L^2[0,l]}\\
\le C_2\|x-x_0\|_{H^1_0(0,l)}\|x\|_{H^1_0(0,l)}+\|\Phi^{\prime}\circ x_0\|_{C[0,l]}\|x-x_0\|_{H^1_0(0,l)}\\
\le (C_2(\|x_0\|_{H^1_0(0,l)}+\varrho)+\|\Phi^{\prime}\circ x_0\|_{C[0,l]})\|x-x_0\|_{H^1_0(0,l)},
\endgathered
\end{equation*}
which together with \eqref{loclip} proves that $f$ is locally Lipschitz.

By means of Theorem 3.3.4 Exercise 1 from the Section 3.3 in \cite{henri}) we can show that the solution to \eqref{ExamplCD} exists globally. To this end it is sufficient to check that
\begin{equation}\label{GlobalEx}
\gathered
\frac{\|f(t,x(\cdot,t),y(t),d_1)\|_{\widehat X}}{1+\|x(t)\|_{H^1_0(0,l)}+\vert y(t)\vert}
\endgathered
\end{equation}
is bounded on the domain of existence of the solution $(x(\cdot,t), y(t))\in X$.

To verify this fact we introduce the following auxiliary function
\begin{equation*}
\gathered
U(x,y)=\|x\|_{L^2[0,l]}^2+y^2,\quad x\in H_0^1[0,l], \quad y\in \Bbb R
\endgathered
\end{equation*}
and use the following
\begin{proposition}\label{prop2}
For any $(x,y,\xi)\in H_0^1(0,l)\times\Bbb R\times U_1$ it holds that
\begin{equation}\label{App1}
\gathered
\dot U(x,y,\xi)\le \zeta^{\T}(\wt A_0+\varepsilon \id)\zeta+\varepsilon^{-1}\|\xi\|^2_{U_1},
\endgathered
\end{equation}
where we denote 
\begin{equation*}
\gathered
\zeta=(\|x\|_{L^2[0,l]},|y|)^{\T},\\
\wt A_0=\begin{pmatrix}
-\frac{2\pi^2a^2}{l^2}&\|B+D\|_{L^2[0,l]}\\
\|B+D\|_{L^2[0,l]}&2c^2
\end{pmatrix},
\endgathered
\end{equation*}
\end{proposition}
The proof can be found in the Appendix.
\begin{corollary}\label{cor1}
Let $(x(\cdot,t),y(t))$ be a solution to \eqref{ExamplCD} with initial conditions
 $x(\cdot,t_0)=x_0$, $y(t_0)=y_0$ and defined on $t\in[t_0, t_0+\epsilon_{t_0,x_0,y_0,d_1})$, then
\begin{equation*}
\gathered
\sup\limits_{t\in[t_0, t_0+\epsilon_{t_0,x_0,y_0,d_1})}\|(x,y)\|_{\widehat X}\le C(x_0,y_0,d_1).
\endgathered
\end{equation*}
\end{corollary}
The proof can be found in the Appendix.

From this corollary and the locally Lipschitzness follows that \eqref{GlobalEx} is bounded, hence the solution
$(x(\cdot,t),y(t))$ of \eqref{ExamplCD} exists for all  $t\ge t_0$.
The properties of $\alpha$, $\beta$ guarantee that $g(H_0^1(0,l)\times\Bbb R\times U_2)\subset H_0^1(0,l)\times\Bbb R$,
which demonstrates the well-posedness of the problem \eqref{Exampl}--\eqref{ExamplBC}.

To state the asymptotic stability conditions of the system \eqref{Exampl}--\eqref{ExamplBC} we define the following symmetric matrices
\begin{equation*}
\gathered
\scriptsize
A_0=
\begin{pmatrix}
-\frac{2\pi^2a^2}{l^2}&\|B\|_{H^1_0(0,l)}+\frac{l}{\pi}\|D\|_{L^2[0,l]}\\
\|B\|_{H^1_0(0,l)}+\frac{l}{\pi}\|D\|_{L^2[0,l]}&2c^2
\end{pmatrix},
\endgathered
\end{equation*}
\begin{equation*}
\gathered
B_0=\\
\scriptsize
\begin{pmatrix}
\|\alpha^2\|_{C[0,l]}+\frac{l^2}{\pi^2}\|\alpha\alpha_{zz}\|_{C[0,l]}+\frac{l^2}{\pi^2}\|\gamma\|_{L^2[0,l]}^2&*\\
\|\alpha_z\beta_z+\gamma\delta\|_{L^2[0,l]}\frac{l}{\pi}+\|\alpha\beta_z\|_{L^2[0,l]}&\delta^2+\|\beta\|_{H_0^1(0,l)}^2
\end{pmatrix},\\
\endgathered
\end{equation*}
\begin{equation*}
\gathered
\sigma=\|\alpha^2\|_{C[0,l]}+\frac{l^2}{\pi^2}\|\alpha\alpha_{zz}\|_{C[0,l]}+\frac{l^2}{\pi^2}\|\gamma\|_{L^2[0,l]}^2\\
+2(\|\alpha_z\beta_z+\gamma\delta\|_{L^2[0,l]}\frac{l}{\pi}+\|\alpha\beta_z\|_{L^2[0,l]})+\delta^2+\|\beta\|_{H_0^1(0,l)}^2,\\
\vartheta:=\frac{\pi^2a^2}{l^2}-\|B\|_{H^1_0(0,l)}-\frac{l}{\pi}\|D\|_{L^2[0,l]}-c^2,
\endgathered
\end{equation*}
let  $\varrho_{\max}^{A_0}\in\Bbb R^2$ and $\varrho_{\max}^{B_0}\in\Bbb R^2$ be eigenvectors of $A_0$  and $B_0$, respectively, corresponding to the maximal eigenvalues $\lambda_{\max}(A_0)$ and $\lambda_{\max}(B_0)$, respectively.
\begin{proposition}\label{prop3}
Let the impulsive system \eqref{Exampl} satisfy
\begin{equation*}
\gathered
c^2+\frac{\pi^2a^2}{l^2}>\frac{l}{\pi}\|D\|_{L^2[0,l]}+\|B\|_{H^1_0(0,l)}, \vartheta>0,\\
(\varrho_{\max}^{A_0})^{\T} \diag\{-1,1\}\varrho_{\max}^{A_0}>0
\endgathered
\end{equation*}
and additionally let one of the two following conditions hold

(a) $(\varrho_{\max}^{B_0})^{\T}\diag\{-1,1\}\varrho_{\max}^{B_0}\ge 0$,
and $\theta_1$, $\theta_2$ satisfy
\begin{equation*}
\gathered
\frac{1}{\vartheta}\ln\frac{\sigma}{2}<\theta_1\le\theta_2<-\frac{1}{\lambda_{\max}(A_0)}\ln\lambda_{\max}(B_0).
\endgathered
\end{equation*}

(b) $(\varrho_{\max}^{B_0})^{\T}\diag\{-1,1\}\varrho_{\max}^{B_0}<0$,
and $\theta_1$, $\theta_2$ satisfy
\begin{equation*}
\gathered
\frac{1}{\vartheta}\ln\lambda_{\max}(B_0)<\theta_1\le\theta_2<\frac{1}{\lambda_{\max}(A_0)}\ln\frac{2}{\sigma}.
\endgathered
\end{equation*}
then \eqref{Exampl} is ISS for all $\mathcal E$ satisfying the dwell-time condition $\theta_1\le T_k\le\theta_2$.
\end{proposition}
To show how the last proposition can be applied we consider the following
\subsection{Specific example} Consider \eqref{Exampl} with $a=1$, $l=\pi$, $D(z)=0.05z$, $B(z)=0.05z(\pi-z)$, $c=0.5$,
$\alpha(z)=1$, $\beta(z)=0$, $\gamma(z)=0.05$, $\delta=0.25$, $\Phi(s)=-s^3$. Then we have
\begin{equation*}
\gathered
A_0=\begin{pmatrix}
-2&0.3214875\\
0.3214875&0.5
\end{pmatrix},\\
B_0=\begin{pmatrix}
1.007853982&0.02215567314\\
0.02215567314&0.0625
\end{pmatrix},
\endgathered
\end{equation*}
\begin{equation*}
\gathered
\lambda_{\max}(A_0)=0.54067976,\quad\lambda_{\max}(B_0)=1.0083729,\\
\sigma=1.1146653,\quad\vartheta=0.42851243,\\
\varrho_{\max}^{A_0}=(0.12553504,-0.9920891862)^{\T},\\
\varrho_{\max}^{B_0}=(0.99972578,-0.023417096)^{\T}
\endgathered
\end{equation*}
Then all conditions of the Proposition \ref{prop3} are satisfied and the dwell-time condition reads as
\begin{equation*}
\gathered
0.01945822<\theta_1\le T_k\le\theta_2<1.0812185.
\endgathered
\end{equation*}
Let us note that with this parameters choice we have that both discrete and continuous dynamics of  \eqref{Exampl} considered separately are not asymptotically stable already for the unperturbed case $d_1=0$, $d_2=0$. Indeed,  $W(x,y)$ can be used as a Chetaev function for the continuous dynamics of \eqref{Exampl}, and to see that the discrete dynamics is unstable just check that the spectral radius of the jump operator is larger than $1$.

\section{Conclusions} Our results provide a dwell-time condition that guarantees the ISS property of a nonlinear impulsive system. In contrary to the existing dwell-time conditions in the literature our result can be in particular applied
even to the cases where both discrete and continuous dynamics are unstable simultaneously. In contrary to the results of  \cite{dash-miron13} the ISS property is assured by the analysis of specific interaction of the discrete and continuous dynamics instead of a compensation of the unstable discrete (continuous) dynamics by means of the stable continuous (discrete) one. Our future research will be devoted to the development of constructive approach in order to derive the auxiliary Lyapunov $V$ and Chetaev $W$ functions. Another open problem that needs to be investigated is the derivation of conditions under which a combination of the simultaneously stable discrete and continuous dynamics leads to an unstable dynamics of the overall impulsive system.

\section{Appendix}
\subsection{Proof of Proposition \ref{prop2}}
Consider the function
\begin{equation*}
\gathered
U(x,y)=\|x\|_{L^2[0,l]}^2+y^2,\quad
\endgathered
\end{equation*}

Its time derivative  $\dot U$ with respect to the system \eqref{Exampl} is
\begin{equation*}
\gathered
\dot U(x,y,\xi)=2\int\limits_0^lx(z)(a^2\partial_{zz} x(z)+\Phi(x(z))+B(z)y\\
+\xi_1(z))\,dz+2y(c^2 y+\int\limits_0^lD(z)x(z)\,dz+\xi_2)
\endgathered
\end{equation*}
Applying integration by parts and the Friedrich's inequality \eqref{Fr1} we get
\begin{equation*}
\gathered
\int\limits_0^lx(z)\partial_{zz}x(z)\,dz=-\int\limits_0^l|\partial_z x(z)|^2\,dz\le-\frac{\pi^2}{l^2}\|x\|_{L^2[0,l]}^2.
\endgathered
\end{equation*}
By means of the Cauchy inequality we can write
\begin{equation*}
\gathered
\dot U(x,y,\xi)\le -\frac{2\pi^2a^2}{l^2}\|x\|_{L^2[0,l]}^2
+2c^2 y^2\\
+2\|B+D\|_{L^2[0,l]}\|x\|_{L^2[0,l]}|y|\\
+2|y||\xi_2|+2\|\xi_1\|_{L^2[0,l]}\|x\|_{L^2[0,l]}\\
\le
\Big(-\frac{2\pi^2a^2}{l^2}+\varepsilon\Big)\|x\|_{L^2[0,l]}^2
+(2c^2+\varepsilon)y^2\\
+2\|B+D\|_{L^2[0,l]}\|x\|_{L^2[0,l]}|y|+
\varepsilon^{-1}(\|\xi_1\|^2_{L^2[0,l]}+\xi_2^2).
\endgathered
\end{equation*}
for all $(x,y,\xi)\in L^2[0,l]\times\Bbb R\times U_1$ and $\varepsilon>0$.

Recall that $\zeta=(\|x\|_{L^2[0,l]},|y|)^{\T}$ hence the last estimate can be written as
\begin{equation}\label{App1}
\gathered
\dot U(x,y,\xi)\le \zeta^{\T}(\wt A_0+\varepsilon \id)\zeta+\varepsilon^{-1}\|\xi\|^2_{U_1},
\endgathered
\end{equation}
which finishes the proof.

\subsection{Proof of the Corollary \ref{cor1}}
Denote  $u(t):=U(x(\cdot,t),y(t))$, $t\in[t_0,t_0+\epsilon_{t_0,x_0,y_0,d_1})$, $d_1\in\mathcal U_1$,
then Proposition \ref{prop2} implies that
\begin{equation*}
\gathered
\dot u(t)\le(\lambda_{\max}(\wt A_0)+\varepsilon)u(t)+\varepsilon^{-1}\|d_1(t)\|_{U_1}^2\\
\le (\lambda_{\max}(\wt A_0)+\varepsilon)u(t)+\varepsilon^{-1}\|d_1\|_{\mathcal U_1}^2.
\endgathered
\end{equation*}
By the comparison principle we obtain
\begin{equation*}
\gathered
u(t)\le e^{(\lambda_{\max}(\wt A_0)+\varepsilon)(t-t_0)}u(t_0)\\
+\varepsilon^{-1}\int\limits_{t_0}^te^{(\lambda_{\max}(\wt A_0)+\varepsilon)(t-s)}\,ds\|d\|_{\mathcal U_1}^2.
\endgathered
\end{equation*}
By the basic inequality $\sqrt{a+b}\le\sqrt{a}+\sqrt{b}$, $a$, $b\ge 0$ we derive
\begin{equation*}
\gathered
\|(x(\cdot,t),y(t))\|_{\widehat X}\le e^{(\lambda_{\max}(\wt A_0)+\varepsilon)\epsilon_{t_0,x_0,y_0,d_1}/2}\|(x_0,y_0)\|_{\widehat X}\\
+
\sqrt{\frac{e^{(\lambda_{\max}(\wt A_0)+\varepsilon)\epsilon_{t_0,x_0,y_0,d_1}}-1}{\varepsilon(\lambda_{\max}(\wt A_0)+\varepsilon)}}\|d_1\|_{\mathcal U_1},
\endgathered
\end{equation*}
for all $t\in[t_0,t_0+\epsilon_{t_0,x_0,y_0,d_1})$ which finishes the proof.

\subsection{Proof of Proposition \ref{prop3}}
Let us define
\begin{equation*}
\gathered
V(x,y)=\|x\|_{H^1_0(0,l)}^2+y^2,\quad
W(x,y)=y^2-\|x\|_{H^1_0(0,l)}^2
\endgathered
\end{equation*}
and calculate
\begin{equation*}
\gathered
\dot V(x,y,\xi)=-2\int\limits_0^l\partial_{zz}x(z)(a^2\partial_{zz} x(z)+\Phi(x(z))\\
+B(z)y+\xi_1(z))\,dz+2y(c^2 y+\int\limits_0^lD(z)x(z)\,dz+\xi_2)
\endgathered
\end{equation*}
Using the integration by parts
\begin{equation}\label{by-parts}
\gathered
\int\limits_0^l\partial_{zz}x(z)\Phi(x(z))\,dz=-\int\limits_0^l\Phi^{\prime}(x(z))\,(\partial_z x(z))^2\,dz,
\endgathered
\end{equation}
\begin{equation*}
\gathered
\int\limits_0^l\partial_{zz}x(z)B(z)\,dz=-\int\limits_0^l\partial_z x(z)\partial_z B(z)\,dz,
\endgathered
\end{equation*}
\begin{equation*}
\gathered
\int\limits_0^l\partial_{zz}x(z)\xi_1(z)\,dz=-\int\limits_0^l\partial_z x(z)\partial_z \xi_1(z)\,dz,
\endgathered
\end{equation*}
as well as the property $\Phi^{\prime}(s)\le 0$ for all $s\in\Bbb R$ and \eqref{Fr2}, we obtain
\begin{equation*}
\gathered
\dot V(x,y,\xi)\le-\frac{2\pi^2a^2}{l^2}\|x\|_{H^1_0(0,l)}^2\\
+2\int\limits_0^l\partial_z x(z)\partial_z B(z)\,dzy
+2\int\limits_0^l\partial_z x(z)\partial_z \xi_1(z)\,dz\\
+2c^2y^2+2\int\limits_0^lD(z)x(z)\,dz y+2\xi_2y
\endgathered
\end{equation*}
By the Friedrich's inequality \eqref{Fr1} we get
\begin{equation}\label{metka}
\gathered
\dot V(x,y,\xi)\le-\frac{2\pi^2a^2}{l^2}\|x\|_{H^1_0(0,l)}^2\\
+2(\|B\|_{H^1_0(0,l)}+\frac{l}{\pi}\|D\|_{L^2[0,l]})\|x\|_{H^1_0(0,l)}|y|\\
+2c^2y^2+2|\xi_2||y|+2\|x\|_{H^1_0(0,l)}\|\xi_1\|_{H^1_0(0,l)}\\
\le\Big(-\frac{2\pi^2a^2}{l^2}+\varepsilon\Big)\|x\|_{H^1_0(0,l)}^2+2(\|B\|_{H^1_0(0,l)}\\
+\frac{l}{\pi}\|D\|_{L^2[0,l]})\|x\|_{H^1_0(0,l)}|y|
+(2c^2+\varepsilon)y^2\\
+\varepsilon^{-1}(\|\xi_1\|_{H^1_0(0,l)}^2+\xi_2^2).
\endgathered
\end{equation}
and by the Cauchy inequality with  \eqref{Fr2} and \eqref{by-parts} we obtain
\begin{equation*}
\gathered
\dot W(x,y,\xi)=2y(c^2y+\int\limits_0^lD(z)x(z)\,dz+\xi_2)\\
+2\int\limits_0^l\partial_{zz}x(z)(a^2\partial_{zz}x(z)+\Phi(x(z))+B(z)y+\xi_1(z))\,dz\\
\ge 2c^2y^2+2a^2\int\limits_0^l(\partial_{zz}x(z))^2\,dz+
2y\int\limits_0^lD(z)x(z)\,dz\\
-2y\int\limits_0^l\partial_z B(z)\partial_z x(z)\,dz
+2y\xi_2-2\int\limits_0^l\partial_zx(z)\partial_z\xi_1(z)\,dz\\
\endgathered
\end{equation*}
Applying \eqref{Fr2} we have
\begin{equation*}
\gathered
\dot W(x,y,\xi)\ge 2c^2y^2+\frac{2\pi^2a^2}{l^2}\|x\|_{H^1_0(0,l)}^2\\
-2(\frac{l}{\pi}\|D\|_{L^2[0,l]}+\|B\|_{H^1_0(0,l)})\|x\|_{H^1_0(0,l)}|y|\\
-2|y||\xi_2|-2\|x\|_{H^1_0(0,l)}\|\xi_1\|_{H^1_0(0,l)}\\
\ge
(2c^2-\varepsilon)y^2+\Big(\frac{2\pi^2a^2}{l^2}-\varepsilon\Big)\|x\|_{H^1_0(0,l)}^2\\
-2(\frac{l}{\pi}\|D\|_{L^2[0,l]}+\|B\|_{H^1_0(0,l)})\|x\|_{H^1_0(0,l)}|y|\\
-\varepsilon^{-1}(\|\xi_1\|^2_{H^1_0(0,l)}+\xi_2^2),
\endgathered
\end{equation*}
for all $(x,y,\xi)\in X\times U_1$. Further we have
\begin{equation*}
\gathered
V(g(x,y,\mu))\le \int\limits_0^l(\alpha_z(z)x(z)+\alpha(z)x_z(z)+\beta_z(z)y\\+\mu_{1z}(z))^2\,dz+
\Big(\int\limits_{0}^l\gamma(z)x(z)\,dz+\delta y+\mu_2\Big)^2\\
\le\|\alpha^2\|_{C[0,l]}\|x_z\|_{L^2[0,l]}^2+\|\beta_z\|_{L^2[0,l]}^2y^2+\|\mu_{1z}\|_{L^2[0,l]}^2\\+
\|\alpha\alpha_{zz}\|_{C[0,l]}\|x\|_{L^2[0,l]}^2+2\|\alpha_z\beta_z+\gamma\delta\|_{L^2[0,l]}\|x\|_{L^2[0,l]}|y|\\
+2\|\alpha\beta_z\|_{L^2[0,l]}\|x\|_{H_0^1(0,l)}|y|
+2\|\alpha_z\|_{C[0,l]}\|\mu_1\|_{H_0^1(0,l)}\|x\|_{L^2[0,l]}\\
+2\|\beta\|_{H^1_0(0,l)}\|\mu_1\|_{H_0^1(0,l)}|y|
+2\|\alpha\|_{C[0,l]}\|\mu_1\|_{H_0^1(0,l)}\|x\|_{H_0^1(0,l)}\\
+\|\gamma\|_{L^2[0,l]}^2\|x\|_{L^2[0,l]}^2+\delta^2y^2+\mu_2^2
+2|\delta||\mu_2||y|\\
+2|\mu_2|\|\gamma\|_{L^2[0,l]}\|x\|_{L^2[0,l]}
\endgathered
\end{equation*}
By the Friedrich's inequality \eqref{Fr1} we have
\begin{equation*}
\gathered
V(g(x,y,\mu))\le
\Big(\|\alpha^2\|_{C[0,l]}+\frac{l^2}{\pi^2}\|\alpha\alpha_{zz}\|_{C[0,l]}+\frac{l^2}{\pi^2}\|\gamma\|_{L^2[0,l]}^2\\
+\varepsilon(\frac{l^2}{\pi^2}\|\alpha_z\|_{C[0,l]}+\|\alpha\|_{C[0,l]}+\|\gamma\|_{L^2[0,l]}\frac{l^2}{\pi^2})\Big)\|x\|_{H^1_0(0,l)}^2\\
+2\Big(\|\alpha_z\beta_z+\gamma\delta\|_{L^2[0,l]}\frac{l}{\pi}+\|\alpha\beta_z\|_{L^2[0,l]}\Big)\|x\|_{H_0^1(0,l)}|y|\\
+(\|\beta\|_{H^1_0(0,l)}^2+\delta^2+\varepsilon(\|\beta\|_{H_0^1(0,l)}+|\delta|))y^2\\
+(1+\varepsilon^{-1}(\|\alpha_z\|_{C[0,l]}+\|\beta\|_{H_0^1(0,l)}+\|\alpha\|_{C[0,l]}))\|\mu_1\|_{H_0^1(0,l)}^2\\
+(1+\varepsilon^{-1}(|\delta|+\|\gamma\|_{L^2[0,l]}))\mu_2^2
\endgathered
\end{equation*}
Let us denote
\begin{equation*}
\gathered
\scriptsize
B_1=
\begin{pmatrix}
\frac{l^2}{\pi^2}\|\alpha_z\|_{C[0,l]}+\|\alpha\|_{C[0,l]}+\|\gamma\|_{L^2[0,l]}\frac{l^2}{\pi^2}&0\\
0&\|\beta\|_{H_0^1(0,l)}+|\delta|
\end{pmatrix},\\
\kappa(\varepsilon)=\max\{1+\varepsilon^{-1}(\|\alpha_z\|_{C[0,l]}+\|\beta\|_{H_0^1(0,l)}+\|\alpha\|_{C[0,l]}),\\
1+\varepsilon^{-1}(|\delta|+\|\gamma\|_{L^2[0,l]})\},
\endgathered
\end{equation*}
then the last inequality can be written as
\begin{equation}\label{App2}
\gathered
V(g(x,y,\mu))\le \zeta^T(B_0+\varepsilon B_1)\zeta+\kappa(\varepsilon)\|\mu\|^2.
\endgathered
\end{equation}

The inequalities \eqref{metka} and \eqref{App2} imply the estimates \eqref{1.2} and \eqref{1.3}:
Denote $w(t)=\|x(\cdot,t)\|_{H^1_0(0,l)}^2+y^2(t)$, $t\ge 0$, then from \eqref{1.2} follows:
\begin{equation*}
\gathered
\dot w(t)\le(\lambda_{\max}(A_0)+\varepsilon)w(t)+\varepsilon^{-1}\|d_1(t)\|_{U_1}^2\\
\le (\lambda_{\max}(A_0)+\varepsilon)w(t)+\varepsilon^{-1}\|d\|_{\mathcal U_1}^2.
\endgathered
\end{equation*}
By the comparison principle we have
\begin{equation*}
\gathered
w(t)\le e^{(\lambda_{\max}(A_0)+\varepsilon)t}w(0)+\\
\varepsilon^{-1}\int\limits_0^te^{(\lambda_{\max}(A_0)+\varepsilon)(t-s)}\,ds\|d_1\|_{\mathcal U_1}^2.
\endgathered
\end{equation*}
and with help of  $\sqrt{a+b}\le\sqrt{a}+\sqrt{b}$, $a$, $b\ge 0$ we get
\begin{equation*}
\gathered
\|(x(\cdot,t),y(t))\|_X\le e^{(\lambda_{\max}(A_0)+\varepsilon)\tau/2}\|(x_0,y_0)\|_X+\\
\sqrt{\frac{e^{(\lambda_{\max}(A_0)+\varepsilon)\tau}-1}{\varepsilon(\lambda_{\max}(A_0)+\varepsilon)}}\|d_1\|_{\mathcal U_1},
\endgathered
\end{equation*}
for all $t\in[0,\tau]$, which proves \eqref{1.2}. Similarly \eqref{1.3} follows from \eqref{App2}.

Let us estimate $\dot V(x,y,\xi)$, $V(g(x,y,\mu))$ on the sets
\begin{equation*}
\gathered
G^+=\{(x,y)\in X\,\,:\,\,|y|\ge \|x\|_{H^1_0(0,l)}\},\\
G^-=\{(x,y)\in X\,\,:\,\,|y|\le \|x\|_{H^1_0(0,l)}\}
\endgathered
\end{equation*}
and $\dot W(x,y,\xi)$ on the set, where $W(x,y)=0$ holds.

By the conditions of Proposition \ref{prop3}  we have\\ $(\varrho_{\max}^{A_0})^{\T}\diag\{-1,1\}\varrho_{\max}^{A_0}>0$ and
\begin{equation*}
\gathered
-\vartheta=-\frac{\pi^2a^2}{l^2}+\|B\|_{H^1_0(0,l)}+\frac{l}{\pi}\|D\|_{L^2[0,l]}+c^2<0.
\endgathered
\end{equation*}
Hence by \eqref{metka} for all $(x,y)\in G^+$ we have
\begin{equation*}
\gathered
\dot V(x,y,\xi)\le (\lambda_{\max}(A_0)+\varepsilon)(\|x\|_{H^1_0(0,l)}^2+y^2)+\varepsilon^{-1}\|\xi\|^2_{U_1}\\
\le(\lambda_{\max}(A_0)+2\varepsilon)(\|x\|_{H^1_0(0,l)}^2+y^2),
\endgathered
\end{equation*}
if $\sqrt{\|x\|_{H^1_0(0,l)}^2+y^2}>\varepsilon^{-1}\|\xi\|_{U_1}$.

For all $(x,y)\in G^-$ we have
\begin{equation*}
\gathered
\dot V(x,y,\xi)\le (-\vartheta+\frac{\varepsilon}{2})(\|x\|_{H^1_0(0,l)}^2+y^2)+\varepsilon^{-1}\|\xi\|^2_{U_1}\\
\le(-\vartheta+\varepsilon)(\|x\|_{H^1_0(0,l)}^2+y^2),
\endgathered
\end{equation*}
if $\sqrt{\|x\|_{H^1_0(0,l)}^2+y^2}>\sqrt{2}\varepsilon^{-1}\|\xi\|_{U_1}$.

Now choosing $\varepsilon>0$ small enough ($\varepsilon<\vartheta$), we can take
$\vph_1(s)=(\vartheta-\varepsilon)s$, $\vph_2(s)=(\lambda_{\max}(A_0)+2\varepsilon)s$.

If the condition (a) is satisfied, that is\\ $(\varrho_{\max}^{B_0})^{\T}\diag\{-1,1\}\varrho_{\max}^{B_0}\ge 0$, then for all $(x,y)\in G^+$ we have
\begin{equation*}
\gathered
V(g(x,y,\mu))\le (\lambda_{\max}(B_0)+\varepsilon\|B_1\|)(\|x\|_{H^1_0(0,l)}^2+y^2)\\
+\kappa(\varepsilon)\|\mu\|^2_{U_2}
\le (\lambda_{\max}(B_0)+\varepsilon(1+\|B_1\|))(\|x\|_{H^1_0(0,l)}^2+y^2).
\endgathered
\end{equation*}
if $\sqrt{\|x\|_{H^1_0(0,l)}^2+y^2}>\sqrt{\varepsilon^{-1}\kappa(\varepsilon)}\|\mu\|_{U_2}$.

For all  $(x,y)\in G^-$ we have
\begin{equation*}
\gathered
V(g(x,y,\mu))\le \Big(\frac{\sigma}{2}+\varepsilon(1+\|B_1\|)\Big)V(x,y),
\endgathered
\end{equation*}
if $\sqrt{\|x\|_{H^1_0(0,l)}^2+y^2}>\sqrt{\varepsilon^{-1}\kappa(\varepsilon)}\|\mu\|_{U_2}$.

In this case $\psi_1(s)=\Big(\frac{\sigma}{2}+\varepsilon(1+\|B_1\|)\Big)s$,
 $\psi_2(s)=(\lambda_{\max}(B_0)+\varepsilon(1+\|B_1\|))s$.

If the condition (b) is satisfied, that is\\  $(\varrho_{\max}^{B_0})^{\T}\diag\{-1,1\}\varrho_{\max}^{B_0}<0$, then for all
 $(x,y)\in G^+$ we have
\begin{equation*}
\gathered
V(g(x,y,\mu))\le \Big(\frac{\sigma}{2}+\varepsilon(1+\|B_1\|)\Big)V(x,y),
\endgathered
\end{equation*}
if  $\sqrt{\|x\|_{H^1_0(0,l)}^2+y^2}>\sqrt{\varepsilon^{-1}\kappa(\varepsilon)}\|\mu\|_{U_2}$.

For all $(x,y)\in G^-$ we have
\begin{equation*}
\gathered
V(g(x,y,\mu))\le (\lambda_{\max}(B_0)+\varepsilon\|B_1\|)(\|x\|_{H^1_0(0,l)}^2+y^2)\\
+\kappa(\varepsilon)\|\mu\|^2_{U_2}
\le (\lambda_{\max}(B_0)+\varepsilon(1+\|B_1\|))(\|x\|_{H^1_0(0,l)}^2+y^2).
\endgathered
\end{equation*}
if $\sqrt{\|x\|_{H^1_0(0,l)}^2+y^2}>\sqrt{\varepsilon^{-1}\kappa(\varepsilon)}\|\mu\|_{U_2}$. In this case
 $\psi_1(s)=(\lambda_{\max}(B_0)+\varepsilon(1+\|B_1\|))s$, $\psi_2(s)=\Big(\frac{\sigma}{2}+\varepsilon(1+\|B_1\|)\Big)s$.

If
\begin{equation*}
\gathered
c^2+\frac{\pi^2a^2}{l^2}>\frac{l}{\pi}\|D\|_{L^2[0,l]}+\|B\|_{H^1_0(0,l)},
\endgathered
\end{equation*}
then $W(x,y)=0$ implies
\begin{equation*}
\gathered
\dot W(x,y,\mu)\ge \Big(c^2+\frac{\pi^2a^2}{l^2}-\big(\frac{l}{\pi}\|D\|_{L^2[0,l]}\\
+\|B\|_{H^1_0(0,l)}\big)-2\varepsilon\Big)(\|x\|_{H_0^1(0,l)}^2+y^2)>0,
\endgathered
\end{equation*}
if $\sqrt{\|x\|_{H_0^1(0,l)}^2+y^2}>\varepsilon^{-1}\|\xi\|_{U_1}$, with small enough $\varepsilon$.

Now applying Theorem \ref{th1} and choosing positive constant $\varepsilon$ and $\delta$ small enough, we see that the conditions of Proposition \ref{prop3} guarantee the  ISS property for the system  \eqref{Exampl}, which finishes the proof.

	\bibliography{impulsive-biblio}
\bibliographystyle{abbrv}

\end{document}